\documentclass[aip,jcp,amsmath,amsfonts,amssymb,preprint]{revtex4-2}

\usepackage{graphicx}
\usepackage{amsthm}
\usepackage{algorithm,algpseudocode}
\usepackage{subcaption}
\usepackage{hyperref}

\newcommand{\etal}{\textit{et~al.\ }}

\newcommand{\matlab}{M\textsc{atlab}}

\newcommand{\order}[1]{\mathcal{O}{(#1)}}
\newcommand{\subindex}[2]{{\begin{subarray}{l} {#1} \\ {#2}\end{subarray}}}

\renewcommand{\c}{\mathrm{c}}
\newcommand{\erf}{\mathrm{erf}}
\newcommand{\erfc}{\mathrm{erfc}}
\newcommand{\erms}{\varepsilon_\mathrm{rms}}
\newcommand{\ii}{\mathrm{i}}
\newcommand{\kintinf}{\bar{k}_\infty}
\newcommand{\midx}{\mathtt{m}}
\newcommand{\nidx}{\mathtt{n}}
\newcommand{\mn}{{\midx\nidx}}
\newcommand{\mnp}{{\mn\v p}}
\newcommand{\nl}{\mathsf{n}}
\renewcommand{\P}{{\mathcal{P}}}
\renewcommand{\S}[3]{\overset{#2}{\underset{#3}{#1}}}
\newcommand{\subR}{_R}
\renewcommand{\v}[1]{\mathbf{#1}}
\newcommand{\vk}{\boldsymbol{k}}
\newcommand{\vka}{\boldsymbol{\kappa}}
\newcommand{\w}[1]{\mathrm{#1}}
\newcommand{\W}{{\mathsf{w}}}
\newcommand{\wh}{\widehat}
\newcommand{\wt}{\widetilde}

\theoremstyle{definition}
\newtheorem{definition}{Definition}

\algnewcommand\algorithmicinput{\textbf{Input:}}
\algnewcommand\INPUT{\item[\algorithmicinput]}
\algnewcommand\algorithmicoutput{\textbf{Output:}}
\algnewcommand\OUTPUT{\item[\algorithmicoutput]}

\makeatletter
\newenvironment{breakablealgorithm}
  {
   \begin{center}
     \refstepcounter{algorithm}
     \hrule height.8pt depth0pt \kern2pt
     \renewcommand{\caption}[2][\relax]{
       {\raggedright\textbf{\ALG@name~\thealgorithm} ##2\par}%
       \ifx\relax##1\relax 
         \addcontentsline{loa}{algorithm}{\protect\numberline{\thealgorithm}##2}%
       \else 
         \addcontentsline{loa}{algorithm}{\protect\numberline{\thealgorithm}##1}%
       \fi
       \kern2pt\hrule\kern2pt
     }
  }{
     \kern2pt\hrule\relax
   \end{center}
  }
\makeatother

\begin{document}

\title{Fast Ewald summation for electrostatic potentials with arbitrary periodicity}

\author{D. S. Shamshirgar}
\email[]{davoudss@kth.se}

\author{J. Bagge}
\email[]{joarb@kth.se}

\author{A.-K. Tornberg}
\email[]{akto@kth.se}
\affiliation{KTH Mathematics, Swedish e-Science Research Centre, 100 44 Stockholm, Sweden.}


\begin{abstract}
  A unified treatment for fast and spectrally accurate evaluation
  of electrostatic potentials subject to periodic boundary
  conditions in any or none of the three spatial dimensions is
  presented. Ewald decomposition is used to split the problem
  into a real-space and a Fourier-space part, and the FFT-based
  Spectral Ewald (SE) method is used to accelerate the
  computation of the latter. A key component in the unified
  treatment is an FFT-based solution technique for the free-space
  Poisson problem in three, two or one dimensions, depending on
  the number of non-periodic directions. The computational cost
  is furthermore reduced by employing an adaptive FFT for the
  doubly and singly periodic cases, allowing for different local
  upsampling factors. The SE method will always be most efficient
  for the triply periodic case as the cost of computing FFTs
  will then be the smallest, whereas the computational cost of
  the rest of the algorithm is essentially independent of
  periodicity. We show that the cost of removing periodic
  boundary conditions from one or two directions out of three
  will only moderately increase the total runtime. Our
  comparisons also show that the computational cost of the SE
  method in the free-space case is around four times that of the
  triply periodic case.

  The Gaussian window function previously used in the SE method,
  is here compared to a piecewise polynomial approximation of the
  Kaiser-Bessel window function. With a carefully tuned shape
  parameter that is selected based on an error estimate for this
  new window function, runtimes for the SE method can be further
  reduced. Furthermore, we consider different methods for
  computing the force, and compare the runtime of the SE method
  with that of the Fast Multipole Method.
\end{abstract}


\maketitle

\section{Introduction}
\label{sec:intro}
The task of computing interactions in an $N$-body problem is the
most demanding part of various numerical simulations such as electrostatics in molecular dynamics, gravitational fields in cosmological formation of galaxies, and potentials in Stokes flow simulations. 
Due to the long-range behavior of the involved kernels, these problems are computationally expensive and therefore, fast and accurate numerical algorithms are required to accelerate simulations.
The Ewald technique \cite{Ewald1921} splits the interactions into a near field (computed in real space) and a far field (computed in Fourier space) contribution. There exist several methods that utilize this decomposition together with the Fast Fourier transform (FFT) in order to accelerate the calculation of the Fourier space sum.\cite{Hockney2010,Deserno1998,Essmann1995,Lindbo2011a}
These methods belong to a family of Particle-Mesh-Ewald (PME) methods which, applied to a system of $N$ particles, reduce the computational complexity from $\order{N^2}$ to $\order{N\log(N)}$ with a prefactor depending on the required accuracy. 

The given references are all concerned with the triply periodic
case. Different approaches have been suggested to extend to settings
with reduced periodicity. One simple idea for the doubly periodic case
is to use the triply periodic summation, simply extending the unit
cell in the non-periodic direction and thereby creating a gap that
separates sheets of charged particles. To increase the accuracy,
various methods have been proposed that introduce correction terms to
the triply periodic sum. In Arnold \etal \cite{Arnold2002b}, a correction term that allows for highly accurate calculations
is derived.

Lekner summation is an alternative to Ewald summation, and the MMM2D
method is based on this approach,\cite{Arnold2002} reducing the computational
complexity not to the desired $O(N \log(N))$ but rather $O(N^{5/3})$.  There
is also a MMM1D method for singly periodic problems,\cite{Arnold2005} but with a computational cost of $O(N^2)$.

In Ref.~\onlinecite{Nestler2015}, $O(N \log(N))$ methods are introduced both for
doubly and singly periodic problems. As soon as a non-periodic
direction exists, the discrete summation in Fourier space in the Ewald
decomposition is substituted with an integral.  Evaluating the
integrals analytically, the summation in Fourier space involves
complementary error functions and Bessel functions for
the doubly and singly periodic cases, respectively. The methods
introduced in Ref.~\onlinecite{Nestler2015}
are based on regularization and periodic extension of such
functions to enable the use of FFTs.

In this paper, we present a different approach, that works directly
with the numerical discretization of the integrals in Fourier space.
The Spectral Ewald method has been developed over the last decade 
\cite{Lindbo2011a,Lindbo2012, Klinteberg2017, Shamshirgar2017} in order to provide a 
fast and spectrally accurate approach
for evaluating electrostatics problems with different periodicities.
The free-space and 1d-periodic versions of the method
were developed recently and equipped with a novel technique proposed
by Vico et al. which provides a tool for computing volume potentials
using FFTs.\cite{Vico2016} Together with an adaptive FFT that enables different
local upsampling factors, this extension makes it possible to unify the
treatment of all modes, as was done in the singly periodic case in Ref.~\onlinecite{Shamshirgar2017}.
As a result, this case can be evaluated with only a small extra cost as 
compared to the triply periodic case.\cite{Shamshirgar2017} In
Ref.~\onlinecite{Klinteberg2017}, the free-space version of the
SE method is used for accelerating the evaluation of free-space
potentials of Stokes flow. It was shown that this method is
competitive with the Fast Multipole Method (FMM), which unlike
the SE method is most efficient for tackling problems with
non-periodic boundary conditions.

In the current work, we extend the recent advances made also to the doubly
periodic electrostatic problem, previously considered in Ref.~\onlinecite{Lindbo2012}, as well
as the free-space case, to complete the full range from free space to
triply periodic in one unified treatment.

The SE method has so far been using a Gaussian window function to
interpolate between point sources and a uniform mesh. In this
paper, we replace the Gaussian window by a piecewise polynomial
approximation of the Kaiser-Bessel (KB) window
function\cite{Kaiser1980} to perform the interpolation. This
approximation, inspired by Ref.~\onlinecite{Barnett2019}, is accurate enough to retain desired properties of
the KB window and is substantially cheaper to evaluate. For both the Gaussian and KB
window function, a shape parameter has to be set. We provide an
error estimate useful for finding the optimal shape parameter of
the KB window function and assess the accuracy of the estimate by
means of numerical tests. We show that employing the polynomial
approximation of the KB window instead of the Gaussian window,
the cost of evaluation using the SE method is reduced
significantly. We also provide a systematic approach for
selecting parameters based on a given error tolerance, which can
be used to automate the parameter selection process.

This paper is organized as follows: In
section~\ref{sec:ewald:sum}, we provide
Ewald summation formulas for different types of periodic boundary
conditions. In section~\ref{sec:SE}, the Spectral Ewald method is
constructed for problems with arbitrary periodicity.
Section~\ref{sec:windows} introduces different window functions
that can be used in PME methods, as well as the piecewise
polynomial approximation used in this paper. Truncation and
approximation errors together with error estimates are introduced
in section~\ref{sec:errors}, while section~\ref{sec:params} is
dedicated to selection of parameters.
The numerical results in section~\ref{sec:results} have a
threefold focus: (i) comparison of the Gaussian and KB window
functions, (ii) computation of forces, which is important in
applications, and (iii) a new comparison with the FMM.
Finally, conclusions are drawn in section~\ref{sec:conclusions}.

\section{Ewald summation}
\label{sec:ewald:sum}

The classical Ewald sum was developed for fast evaluation of
potentials in ionic crystals and later the same technique was
used for computing long-range interactions in molecular dynamics
simulations and potentials in Stokes flow. The resulting formula
relies on the Ewald decomposition introduced by
Ewald\cite{Ewald1921} in 1921 for 3d-periodic problems. The Ewald
sum in the 2d-periodic case, sometimes referred to as
\textit{slab/slablike geometry}, was derived by Grzybowski \etal
\cite{Grzybowski2000} using lattice sums. The first derivation of
the Ewald sum for the 1d-periodic case was given by Porto
\cite{Porto2000} using an integral representation of the Gamma
function and the Poisson summation formula. The author left an
integral in his expression; however, following
Ref.~\onlinecite{Fripiat2010}, the closed form of the integral
can be obtained. For alternative derivations of Ewald formulas
with different periodicities, the reader may consult
Ref.~\onlinecite{Tornberg2014}.

Consider a system of $N$ point sources with charges $q_\nidx \in
\mathbb{R}$ located at positions $\v x_\nidx \in \Omega$,
$\nidx=1,2,\ldots,N$, in a cubic box $\Omega=[0,L)^3 \subset
\mathbb{R}^3$. The objective is to calculate the electrostatic
potential $\varphi$, given by the discrete sum
\begin{align}
\varphi(\v x_\midx)=\sum_{\v p\in \P_D}^{'}\sum_{\nidx=1}^N
  \frac{q_\nidx}{|\v x_\midx - \v x_\nidx + \v p|},
\qquad \midx=1,2,\ldots,N.
\label{eq:potential}
\end{align}
The prime above the summation
symbol denotes that the term with $\nidx=\midx$ and $\v p=\v 0$ is
excluded from the sum. The set $\P_D$ with
$D\in\lbrace0,1,2,3\rbrace$ is defined to impose periodicity,
\begin{align*}
\begin{array}{ll}
\text{Triply periodic}: & \P_3 = \lbrace(\alpha_1L,\alpha_2L,\alpha_3L): \alpha_i\in\mathbb{Z}\rbrace,\\
\text{Doubly periodic}: & \P_2 = \lbrace(\alpha_1L,\alpha_2L,0): \alpha_i\in\mathbb{Z}\rbrace,\\
\text{Singly periodic}: & \P_1 = \lbrace(\alpha_1L,0,0): \alpha_i\in\mathbb{Z}\rbrace,\\
\text{Free space}: & \P_0 = \lbrace(0,0,0)\rbrace.
\end{array}
\end{align*}
The computational domain~$\Omega$ does not have to be a cube but this
assumption simplifies the description and formulation.
We also assume that the system is charge-neutral, i.e.\
$\sum_\nidx q_\nidx=0$.
This condition is necessary for the sum to converge in
triply, doubly and singly periodic cases\cite{Smith2008};
however, we assume that it also holds for free-space systems.
Even for a charge-neutral system, the sum in
Eq.~\eqref{eq:potential} is only conditionally convergent in the
triply periodic case and therefore the order of summation has to
be defined.\cite{deLeeuw1980} The classical Ewald summation
formula derived in Ref.~\onlinecite{Ewald1921} corresponds to a
spherical order of summation and so-called ``tin foil'' far-field
conditions, i.e.\ a surrounding medium with infinite dielectric
constant.\cite{Tornberg2014}

The potential \eqref{eq:potential} is the solution to the Poisson
problem
\begin{align}
-\Delta\varphi(\v x) = f^{D\P}(\v x), \quad f^{D\P}(\v x) = 4\pi
  \sum_{\v p\in \P_D}\sum_{\nidx=1}^Nq_\nidx \delta(\v x-\v
  x_\nidx+\v p),
  \qquad \v x \in \mathbb{R}^3,
\label{eq:laplace}
\end{align}
with the conditions that $\nabla \varphi(\v x)$ vanishes at
infinity in the free (i.e.\ non-periodic) directions and that
\begin{equation}
  \int_{\mathbb{R}^{3-D}} \int_{[0,L)^D} \varphi(\v x) \,
  \textup{d} \v x = 0,
\end{equation}
where the integral is over $\mathbb{R}$ in each free direction
and over $[0,L)$ in each periodic direction.
In Eq.~\eqref{eq:laplace}, $\Delta$ is the Laplace operator and
$\delta$ is the Dirac delta function.
By introducing a screening function $\gamma$,
$f^{D\P}$ is decomposed as
\begin{align}
f^{D\P} = f^{D\P,\w R}+f^{D\P,\w F}, \quad f^{D\P,\w R} =
  f^{D\P}-(f^{D\P}*\gamma), \quad f^{D\P,\w F} =
  (f^{D\P}*\gamma),
\label{eq:split_rhs}
\end{align}
where $*$ denotes convolution. Now, $\varphi(\v x_\midx)$ can be
obtained by solving two Poisson equations with the right-hand
sides  of $f^{D\P,\w R}$ and $f^{D\P,\w F}$. The solutions to
these two problems are denoted here by $\varphi^{D\P,\w R}(\v
x_\midx)$ and $\varphi^{D\P,\w F}(\v x_\midx) $ respectively.  The
total solution to the problem in Eq.~\eqref{eq:laplace} can then
be written as
\begin{align}
\varphi(\v x_\midx) = \varphi^{D\P,\w R}(\v x_\midx) +
  \varphi^{D\P,\w F}(\v x_\midx) + \varphi^{\w {self}}_\midx.
\label{eq:potential_split}
\end{align}
To obtain the classical Ewald sum, the screening function
$\gamma$, with Fourier transform $\hat{\gamma}$, is selected as
\begin{align}
  \gamma(\v x,\xi)=\xi^3\pi^{-3/2}e^{-\xi^2|\v x|^2}, \qquad 
  \hat{\gamma}(\v k,\xi) = e^{-|\v k|^2/4\xi^2}.
  \label{eq:gamma}
\end{align}
Here, the Ewald decomposition parameter $\xi > 0$ controls how
fast the two terms $\varphi^{D\P,\w R}(\v x_\midx)$ and
$\varphi^{D\P,\w F}(\v x_\midx)$ decay, but does not change the
final result $\varphi(\v x_\midx)$.

The self-contribution term $\varphi^{\w {self}}_\midx$ in
Eq.~\eqref{eq:potential_split} is a constant term which is
independent of the periodicity. This term is added to the sum in
order to exclude the unwanted interaction of charges with
themselves which is introduced as a result of the decomposition. The term reads
\begin{align*}
\varphi^{\w {self}}_\midx = -\dfrac{2\xi}{\sqrt{\pi}}q_\midx.
\end{align*}
The real-space sum can be written as
\begin{align}
\varphi^{D\P,\w R}(\v x_\midx) =  \sum_{\v p\in
  \P_D}^{'}\sum_{\nidx=1}^N  q_\nidx\dfrac{\erfc(\xi|\v
  x_\mnp|)}{|\v x_\mnp|},\quad \midx=1,\ldots,N,
\label{eq:real_space_ewald}
\end{align}
where $\erfc$ is the complementary error function and
\begin{align*}
\v x_\mnp := \v x_\mn+\v p := \v x_\midx-\v x_\nidx+\v p,
\end{align*}
and as before $D\in\lbrace0,1,2,3\rbrace$ represents free-space,
1d-, 2d- and 3d-periodic cases. The sum in
Eq.~\eqref{eq:real_space_ewald} decays exponentially fast with
$|\v x_\mnp|$ and is calculated by introducing a cut-off radius
$r_\c>0$ and including only terms s.t. $|\v x_\mnp|<r_\c$. In
practice, a cell list is constructed for each target point $\v
x_\midx$. Taking into account that the domain is wrapped around periodically in periodic directions, the calculation is restricted to this list. 

The term $\varphi^{D\P,\w F}(\v x_\midx)$ is smooth, and therefore
its Fourier spectrum decays rapidly. This term is treated in
Fourier space, and its structure depends heavily on the type of
periodicity, i.e.\ on $D$. In order to unify the description for
different periodicities, we now introduce some non-standard notation.

\begin{definition}
  Let $D \in \{0, 1, 2, 3\}$ be the number of periodic directions
  and write $\v x=[\v v, \v w] = (x,y,z)$ for the spatial position and $\v k=[\vk,\vka] =
  (k_1,k_2,k_3)$ for the wavenumber vector, where $\v w, \vka\in
  \mathbb{R}^{3-D}$ represent free directions and $\v
  v\in\mathbb{R}^D$, $\vk\in\mathcal{K}^D$ represent periodic
  directions with
  \begin{equation*}
    \mathcal{K}^D := \lbrace \vk \in \mathbb{R}^D :
    k_i\in\frac{2\pi}{L}\mathbb{Z},i=1,\ldots,D\rbrace.
  \end{equation*}
  For $D=3$, $\v w$ and $\vka$ are not defined and $\v k=\vk$,
  $\v x=\v v$. For $D=0$, $\v v$ and $\vk$ are not defined and
  $\v k=\vka$, $\v x=\v w$. Also write $k := \lvert \vk \rvert$ and
  $\kappa := \lvert \vka \rvert$.
  Furthermore, define a functional $\mathcal{L}$ that takes a
  function $g : \mathcal{K}^D \times \mathbb{R}^{3-D} \to
  \mathbb{C}$ by
  \begin{equation}
    \mathcal{L}[g(\v k)]=\mathcal{L}[g([\vk,\vka])] :=
    \left\lbrace
    \begin{array}{cl}
      \displaystyle\dfrac{1}{L^D}\S{\sum}{}{\v k\in \mathcal{K}^D}
      g(\v k), &\quad D=3, \\[7mm]
      \displaystyle\dfrac{1}{(2\pi)^{3-D}L^D}\S{\sum}{}{\vk\in
      \mathcal{K}^D}\S{\int}{}{\mathbb{R}^{3-D}} g([\vk,\vka])
      \, \textup{d} \vka,& \quad D\in\lbrace1,2\rbrace, \\[7mm]
      \displaystyle\dfrac{1}{(2\pi)^{3-D}}\S{\int}{}{\mathbb{R}^{3-D}}
      g(\v k) \, \textup{d} \v k,& \quad D=0.
    \end{array}
    \right.
    \label{eq:operator}
  \end{equation}
  \hfill$\blacklozenge$
  \label{def:definition_mixed}
\end{definition}

Let $f([\v v,\v w])$ be a function that is periodic in $\v v$ and
non-periodic in $\v w$ with Fourier transform given by
\begin{equation}
  \hat{f}([\vk,\vka]) := \int_{\mathbb{R}^{3-D}} \int_{[0,L)^D}
  f([\v v, \v w]) e^{-\ii \vk \cdot \v v} e^{-\ii \vka \cdot \v w}
  \, \textup{d} \v v \textup{d} \v w.
\end{equation}
Then $f$ and $\hat{f}$ are related through
\begin{equation}
f([\v v,\v w]) = \mathcal{L}[\hat{f}([\vk,\vka])e^{\ii \vk\cdot \v v} e^{\ii \vka\cdot \v w}].
\label{eq:inverse_transform}
\end{equation}

Using the notation introduced in
definition~\ref{def:definition_mixed}, the $k$-space part of the
Ewald sum reads
\begin{align}
\varphi^{D\P,\w F}(\v x_{\midx})=4\pi\sum_{\nidx=1}^N
  q_{\nidx}\mathcal{L}\left[\dfrac{e^{-(k^2+\kappa^2)/4\xi^2}}{k^2+\kappa^2}e^{\ii
  \vk\cdot(\v v_\midx-\v v_\nidx)} e^{\ii \vka\cdot(\v w_\midx-\v
  w_\nidx)}\right], \qquad \midx=1,\ldots,N.
  \label{eq:fourier_space_nonzero_ewald}
\end{align}
While the compact notation here will help us to unify the
treatment of all periodicities, we realize that it may be
unfamiliar and therefore write out
Eq.~\eqref{eq:fourier_space_nonzero_ewald} explicitly for all
periodicities at the beginning of section~\ref{sec:SE}, see
Eqs.~\eqref{eq:fourier_space_explicit_3P}--\eqref{eq:fourier_space_explicit_0P}.

For $D=3$, the term corresponding to $\v k=\v 0$ in
Eq.~\eqref{eq:fourier_space_nonzero_ewald} vanishes under the
assumed spherical order of summation, charge neutrality and tin
foil conditions.\cite{Ewald1921} For $D=0$, the operator
$\mathcal{L}$ only includes Fourier integrals defined for all
$\vka\in \mathbb{R}^3$ and the closed form of the integral is
nothing but the complement of the real-space sum minus the self
term.

For $D\in\lbrace1,2\rbrace$, the Fourier integrals in
Eq.~\eqref{eq:fourier_space_nonzero_ewald} are defined
for discrete modes $\vk\in \mathcal{K}^D$. For $\vk \neq \v0$ these integrals can be evaluated analytically. We have
\begin{align}
\varphi^{2\P,\w F,\vk\neq\v0}(\v
  x_\midx)=&\dfrac{\pi}{L^2}\sum_{\nidx=1}^N\sum_{\vk\neq0}q_\nidx\frac{e^{\ii\vk\cdot\v v_\mn}}{k}\left[e^{k z_\mn}\erfc\left(\dfrac{k}{2\xi}+\xi z_\mn\right)
e^{-k z_\mn}\erfc\left(\dfrac{k}{2\xi}-\xi z_\mn\right)\right],\label{eq:2pno0}\\
\varphi^{1\P,\w F,\vk\neq\v0}(\v
  x_\midx)=&\frac{1}{L}\sum_{\nidx=1}^N\sum_{k_1\neq0}q_{\nidx}e^{\ii k_1x_\mn}{\v K_0}(k_1^2/4\xi^2,|\v w_\mn|^2\xi^2),
\label{eq:1pno0}
\end{align}
in which we used the fact that for $D=2$, $\v w_\mn=z_\mn$ and
for $D=1$, $\v v_\mn=x_\mn$ and $\vk=k_1$. Here, ${\v
K_0}(\cdot,\cdot)$ is the incomplete modified Bessel function of
the second kind, defined as
\begin{equation*}
  \v K_0(a,b) = \int_1^\infty \frac{e^{-at-b/t}}{t} \, \w d t.
\end{equation*}
Note that, as done in Ref.~\onlinecite{Nestler2016}, it is
possible to construct a fast method for the sums defined in
Eqs.~\eqref{eq:2pno0} and \eqref{eq:1pno0}. However, using the integral representation of the sums \eqref{eq:fourier_space_nonzero_ewald}, we are able to construct a fast method that has a minimal deviation from the treatment of the triply periodic SE method \cite{Lindbo2011a} while incurring only a small additional cost.

For $D\in\lbrace1,2\rbrace$ and $\vk=\v0$, the Fourier integrals in
Eq.~\eqref{eq:fourier_space_nonzero_ewald} are singular but have closed-form solutions
\begin{align}
\varphi^{2\P,\w F,\vk=\v0}(\v
  x_\midx)=&-\dfrac{2\sqrt{\pi}}{L^2}\sum_{\nidx=1}^Nq_\nidx\left[ e^{-\xi^2z_\mn^2}/\xi+\sqrt{\pi}z_\mn \erf(\xi z_\mn) \right] ,\label{eq:2p0}\\
\varphi^{1\P,\w F,\vk=\v0}(\v
  x_{\midx})=&-\dfrac{1}{L}\sum_\subindex{\nidx=1}{\nidx\neq
  \midx}^Nq_{\nidx}\lbrace \gamma+\log(\xi^2 |\v w_\mn|^2)+\w E_1(\xi^2 |\v w_\mn|^2)\rbrace,
\label{eq:1p0}
\end{align}
where $\erf$ is the error function, $\w E_1(\cdot)=\v K_0(\cdot,0)$ and $\gamma=0.5772156649\ldots$ is the Euler-Mascheroni constant.
The term $\varphi^{2\P,\w F,\vk=\v0}$, cf.\ Eq.~\eqref{eq:2p0}, is a
one-dimensional sum and in Ref.~\onlinecite{Lindbo2012} it is
computed via Chebyshev interpolation. This approach would be much
more expensive if it were to be used for the two-dimensional sum
$\varphi^{1\P,\w F,\vk=\v0}$. In fact, using Chebyshev
interpolation in this case, the cost of the zero-mode term
\eqref{eq:1p0} would be comparable to the cost of calculating the
rest of the Fourier modes, cf. Eq.~\eqref{eq:1pno0}. Instead, we
note that for $\vk=\v0$, \eqref{eq:2p0} and \eqref{eq:1p0} are
solutions to $(3-D)$-dimensional free-space Poisson problems.
Using the idea in Ref.~\onlinecite{Vico2016}, the corresponding
forms in Eq.~\eqref{eq:fourier_space_nonzero_ewald} can be
replaced by non-singular expressions amenable to numerical
integration. This approach has already been used in
Ref.~\onlinecite{Shamshirgar2017,Klinteberg2017} for designing
the Spectral Ewald method for the 1d-periodic and free-space
cases. In this paper, the same technique is used for the zero
mode of the 2d-periodic case as well. In Sec.~\ref{sec:vico} we
review the treatment of the $\vk=\v0$ case in detail.

\section{The Spectral Ewald method}
\label{sec:SE}

Here, we introduce the unified Spectral Ewald (SE) method, a fast
method to accelerate the computation of the Fourier-space part of
the electrostatic potential, i.e.\
Eq.~\eqref{eq:fourier_space_nonzero_ewald}. To aid the reader, we
first write down Eq.~\eqref{eq:fourier_space_nonzero_ewald} explicitly
for the separate cases $D=3,2,1,0$,
\begin{align}
  \varphi^{3\P,\w F}(\v x_{\midx}) &= 4\pi \sum_{\nidx=1}^N
  q_{\nidx} \frac{1}{L^3} \sum_{\substack{(k_1,k_2,k_3) \in \mathcal{K}^3\\[2pt]
  (k_1,k_2,k_3) \neq \v0}}
  \frac{e^{-(k_1^2+k_2^2+k_3^2)/4\xi^2}}{k_1^2+k_2^2+k_3^2}
  e^{\ii (k_1,k_2,k_3) \cdot(\v x_\midx-\v x_\nidx)},
  \label{eq:fourier_space_explicit_3P}\\
  \varphi^{2\P,\w F}(\v x_{\midx}) &= 4\pi \sum_{\nidx=1}^N
  q_{\nidx} \frac{1}{2\pi L^2}
  \sum_{(k_1,k_2) \in \mathcal{K}^2}
  \int_{\mathbb{R}}
  \frac{e^{-(k_1^2+k_2^2+\kappa_3^2)/4\xi^2}}{k_1^2+k_2^2+\kappa_3^2}
  e^{\ii (k_1,k_2,\kappa_3) \cdot (\v x_\midx-\v x_\nidx)} \, \textup{d} \kappa_3,
  \label{eq:fourier_space_explicit_2P}\\[12pt]
  \varphi^{1\P,\w F}(\v x_{\midx}) &= 4\pi \sum_{\nidx=1}^N
  q_{\nidx} \frac{1}{(2\pi)^2 L}
  \sum_{k_1 \in \mathcal{K}^1}
  \int_{\mathbb{R}^2}
  \frac{e^{-(k_1^2+\kappa_2^2+\kappa_3^2)/4\xi^2}}{k_1^2+\kappa_2^2+\kappa_3^2}
  e^{\ii (k_1,\kappa_2,\kappa_3)\cdot(\v x_\midx-\v x_\nidx)}
  \, \textup{d} \kappa_2 \textup{d} \kappa_3,
  \label{eq:fourier_space_explicit_1P}\\[12pt]
  \varphi^{0\P,\w F}(\v x_{\midx}) &= 4\pi \sum_{\nidx=1}^N
  q_{\nidx} \frac{1}{(2\pi)^3}
  \int_{\mathbb{R}^3}
  \frac{e^{-(\kappa_1^2+\kappa_2^2+\kappa_3^2)/4\xi^2}}{\kappa_1^2+\kappa_2^2+\kappa_3^2}
  e^{\ii (\kappa_1,\kappa_2,\kappa_3)\cdot(\v x_\midx-\v x_\nidx)}
  \, \textup{d} \kappa_1 \textup{d} \kappa_2 \textup{d} \kappa_3.
  \label{eq:fourier_space_explicit_0P}
\end{align}
Here, we have taken the liberty to write $\kappa_i$ rather than
$k_i$ for the wavenumbers in the free directions, to emphasize
which directions are free. Note that the kernel $(e^{-\lvert \v k
\rvert^2/4\xi^2}/\lvert \v k\rvert^2) e^{\ii \v k \cdot (\v
x_\midx
- \v x_\nidx)}$ is the same for all periodicities; the difference
is that wavenumbers are summed in the periodic directions but
integrated in the free directions. In the following, we will use the
compact form in Eq.~\eqref{eq:fourier_space_nonzero_ewald} to
treat all periodicities simultaneously.

However, we first point out some difficulties related to the
factor $1/\lvert \v k\rvert^2$. In the SE method,
we will discretize the integrals in
Eqs.~\eqref{eq:fourier_space_explicit_2P}--\eqref{eq:fourier_space_explicit_0P}
using the trapezoidal rule on a uniform grid in $k$-space, so
that we may use the Fast Fourier Transform (FFT) to go between
real space and $k$-space. However, note that when the periodic
wavenumber vector $\vk=(k_1,\ldots,k_D)$ is zero in the $D=1,2$
cases, and always in the $D=0$ case, the integrand becomes
singular at the point $\vka=(k_{D+1},\ldots,k_3)=\v0$. Moreover,
when $\vk$ is non-zero but small in the $D=1,2$ cases, the
integrand will vary rapidly when $\vka$ is close to zero, and
will therefore be hard to resolve. The former problem (singular
integrand) is treated by modifying the Green's function $1/\lvert
\v k \rvert^2$ as in Ref.~\onlinecite{Vico2016} to remove the
singularity, which we describe in section~\ref{sec:vico}. The
latter problem (rapidly varying but non-singular integrand) is
treated by adaptive upsampling, described in
section~\ref{sec:aft}.

The SE method was first introduced by Lindbo and Tornberg
\cite{Lindbo2011a} for the 3d-periodic case and was extended by
the same authors for the 2d-periodic case \cite{Lindbo2012}. The
method follows the same steps as other PME methods, namely:
\begin{enumerate}
  \item
    (Gridding) In real space, a uniform grid is introduced and
    irregular point sources are distributed onto the grid using
    an interpolating \emph{window} function. In the SE method,
    Gaussians have traditionally been used as window functions.
  \item
    (FFT) An FFT is applied to compute the Fourier transform of
    the gridded function.
  \item
    (Scaling) In Fourier space, the result is scaled with a
    (modified) Green's function.
  \item
    (IFFT) An IFFT is employed to take the result back to real
    space.
  \item
    (Gathering) Finally, the Fourier-space part of the potential,
    cf.~\eqref{eq:fourier_space_nonzero_ewald}, is evaluated at
    target points by interpolation with the same window function.
\end{enumerate}
The choice of window function influences the accuracy and runtime
of the resulting method. The main feature of the SE method
compared to other PME methods is that the support of the window
function can be varied independently of the size of the uniform
grid, which allows approximation errors stemming from the window
function to be controlled separately from truncation errors (see
section~\ref{sec:errors}).

In the following, we first introduce the modified Green's
functions to treat the case when the integrand is singular
(section~\ref{sec:vico}). Then, we give an outline of
the SE method and how the window function is introduced in the
method (section~\ref{sec:SE_outline}). This is followed by the
discretization (section~\ref{sec:disc}), adaptive FFT and
upsampling (section~\ref{sec:aft}), precomputation for the
free-space case (section~\ref{sec:precomputation}), and finally a
summary of the SE algorithm (section~\ref{sec:SE_summary}) in
Algorithm~\ref{alg:spectral_ewald}.

\subsection{Modified Green's functions for the singular case}
\label{sec:vico}

As noted above, the Fourier integrals in the free directions
require special treatment when $\vk=\v0$ for $D=1,2$, and also
when $D=0$, since the integrand is singular. In this section, we
describe how the singularity can be removed by modifying the
Green's function. This is done by considering a free-space
Poisson problem. As mentioned in connection to
Eqs.~\eqref{eq:2p0}--\eqref{eq:1p0}, the Fourier integrals
corresponding to $\vk=\v0$ for $D=1,2$ can be seen as solutions
to a $(3-D)$-dimensional free-space Poisson problem. For $D=0$,
the whole solution $\varphi^{0\P,\w F}$ is naturally also a
solution to such a problem and can be treated in the same way.
The main ideas are introduced below. We refer the
reader to the original reference \cite{Vico2016} for more
details, to Ref.~\onlinecite{Klinteberg2017} for an application
in three dimensions ($D=0$) and to
Ref.~\onlinecite{Shamshirgar2017} for an application in two
dimensions ($D=1$). Below, we explain the technique in three
dimensions, i.e.\ for $D=0$.

Consider the free-space Poisson problem
\begin{equation}
  -\Delta \varphi(\v x) = f(\v x), \qquad \v x \in \mathbb{R}^3,
  \label{eq:free_space_poisson}
\end{equation}
with boundary conditions $\varphi(\v x)\to 0$ as $|\v
x|\to\infty$. We are interested in the solution $\varphi$ in a
box $\Omega=[0,L)^3 \subset \mathbb{R}^3$, which corresponds to
the box introduced above Eq.~\eqref{eq:potential}.
The solution to problem~\eqref{eq:free_space_poisson} can be
written in real space as a convolution between the Green's
function $G(\v x)$ and the right hand side $f(\v x)$, or in
Fourier space, as
\begin{align}
  \varphi(\v x) = \int_{\mathbb{R}^3} G(\v x-\v y)f(\v y) \, \w d \v y
  = \dfrac{1}{(2\pi)^3} \int_{\mathbb{R}^3}\wh{G}(\v k)\hat{f}(\v
  k)e^{\ii \v k \cdot \v x} \, \w d \v k,
\label{eq:free_space}
\end{align}
where $G(\v x)=1/(4\pi |\v x|)$ and $\wh{G}(\v k) = 1/|\v k|^2$
is the Fourier transform of $G$. Now assume that $f(\v x)$ is
compactly supported in an extended domain
$\tilde{\Omega}=[-\tfrac{1}{2}\delta_L,L+\tfrac{1}{2}\delta_L)^3$
for some $\delta_L \geq 0$, let $\tilde{L} = L + \delta_L$ and
define the diameter
$R=\w{diam}(\tilde{\Omega})=\sqrt{3}\tilde{L}$.
Noting that $\v x$ is always contained in $\tilde{\Omega}$ in
\eqref{eq:free_space}, we can, without modifying the value of
$\varphi(\v x)$, replace $G(\v x) = G(|\v x|) = G(r)$ with a
truncated version
\begin{align*}
G\subR(r) = G(r) \, \w {rect}\!\left(\dfrac{r}{R}\right),
\end{align*}
where
\begin{align*}
  \text{rect}(x) = \left\lbrace
  \begin{array}{cl}
  1,& |x|\le 1,\\
  0,& |x|>1.\\
  \end{array}
  \right.
\end{align*}
Then, using the fact that $G\subR(r)$ is radially symmetric, its
Fourier transform can be computed as\cite{Vico2016}
\begin{align*}
\wh{G}\subR(k) =2\left(\dfrac{\sin(Rk/2)}{k}\right)^2,\\
\end{align*}
and in the limit $k \to 0$ we have
\begin{align*}
\wh{G}\subR(0)=\lim_{k\to0}\wh{G}\subR(k)=\dfrac{R^2}{2}.
\end{align*}
Note that $\wh{G}\subR$ is not singular.
The expression here is valid for the case $D=0$. Using the same
technique, we can derive corresponding Fourier transforms also
for the cases $D=1$ and $D=2$, where the free-space Poisson
equations are two- and one-dimensional, respectively.

Dropping the subscript $R$ and recalling that $\wh{G}(\v k)=1/|\v
k|^2$ in the non-singular case when $\vk \neq \v0$, we define
modified Green's functions $\wh{G}^{D\P}(\v k)$ as follows, with
$R = \sqrt{3-D}\tilde{L}$:
\begin{itemize}
  \item
    For $D=0$ (free in all directions), $\vka=\v k=(k_1,k_2,k_3)$
    and $\vk$ is not defined, and
    \begin{align}
      \wh{G}^{0\P}(\v k)=\left\lbrace
      \begin{array}{ll}
      2\left(\dfrac{\sin(R|\v k|/2)}{|\v k|}\right)^2, & \v k\neq \v 0, \\
      \vspace{-.3cm}\\
      \dfrac{R^2}{2}, & \v k= \v 0, \\
      \end{array}
      \right.
      \label{eq:G_hat_0d}
    \end{align}
    as we showed above.

  \item
    For $D=1$ (periodic in the $x$ direction and free in the $y$
    and $z$ directions), $\vk=k_1$ and $\vka=(k_2,k_3)$, and%
    \cite{Vico2016,Shamshirgar2017}
    \begin{align}
      \wh{G}^{1\P}(\v k)=\wh{G}^{1\P}([\vk,\vka])=\left\lbrace
      \begin{array}{ll}
      \dfrac{1}{|\v k|^2}, & \vk\neq \v 0, \\
      \vspace{-.3cm}\\
      \dfrac{1-J_0(R|\vka|)}{|\vka|^2}-\dfrac{R\log(R)J_1(R|\vka|)}{|\vka|},
        & \vk=\v 0,~\vka\neq \v 0,\\
      \vspace{-.3cm}\\
      \dfrac{R^2}{4}(1-2\log(R)), & \v k= \v 0, \\
      \end{array}
      \right.
      \label{eq:G_hat_1d}
    \end{align}
    where $J_0$ and $J_1$ are the Bessel functions of the first
    kind and order 0 and 1, respectively.

  \item
    For $D=2$ (periodic in the $x$ and $y$ directions and free in
    the $z$ direction), $\vk=(k_1,k_2)$ and $\vka=k_3$, and (see
    Appendix~\ref{app:G_hat_2d})
    \begin{align}
      \wh{G}^{2\P}(\v k)=\wh{G}^{2\P}([\vk,\vka])=\left\lbrace
      \begin{array}{ll}
      \dfrac{1}{|\v k|^2}, & \vk\neq \v 0, \\
      \vspace{-.3cm}\\
      \dfrac{1 - \cos(R|\vka|) - R |\vka| \sin(R|\vka|)}{|\vka|^2},
        & \vk=\v 0,~\vka\neq \v 0,\\
      \vspace{-.3cm}\\
      -\dfrac{R^2}{2}, & \v k= \v 0. \\
      \end{array}
      \right.
      \label{eq:G_hat_2d}
    \end{align}

  \item
    For $D=3$ (periodic in all directions), $\vk=\v
    k=(k_1,k_2,k_3)$ and $\vka$ is not defined, and
    \begin{align}
      \wh{G}^{3\P}(\v k)=\left\lbrace
      \begin{array}{ll}
      \dfrac{1}{|\v k|^2}, & \v k\neq \v 0, \\
      \vspace{-.3cm}\\
      0, & \v k= \v 0, \\
      \end{array}
      \right.
      \label{eq:G_hat_3d}
    \end{align}
    which is just the standard Green's function modified to take
    into account that the $\v k = \v 0$ term vanishes for $D=3$.
\end{itemize}

In the derivation above, we have assumed the right-hand side to
be compactly supported in the extended box $\tilde{\Omega}$. In
\eqref{eq:split_rhs}, however, the right-hand
side is a superposition of screening functions, i.e. Gaussians,
and hence does not have compact support. Nonetheless, in
practice, it decays rapidly outside of $\Omega=[0,L)^3$ and the
extended domain $\tilde{\Omega}$ can be selected such that the
magnitude is arbitrarily small (we return to how this is done in
section~\ref{sec:disc}). Having this in mind, we may now replace
$1/(k^2+\kappa^2)$ by $\wh{G}^{D\P}$ in
\eqref{eq:fourier_space_nonzero_ewald} and define
\begin{align}
\wt{\varphi}^{D\P,\w F}(\v x_{\midx}):=4\pi\sum_{\nidx=1}^N
  q_{\nidx}\mathcal{L}\left[e^{-(k^2+\kappa^2)/4\xi^2}\wh{G}^{D\P}([\vk,\vka])e^{\ii
  \vk\cdot(\v v_\midx-\v v_\nidx)} e^{\ii \vka\cdot(\v w_\midx-\v
  w_\nidx)}\right].
  \label{eq:fourier_space_ewald}
\end{align}
In the triply periodic case, \eqref{eq:fourier_space_ewald} is
just another representation of
\eqref{eq:fourier_space_nonzero_ewald}, and $\wt{\varphi}^{3\P,\w
F}(\v x_{\midx})=\varphi^{3\P,\w F}(\v x_{\midx})$. In the cases
$D=\{0,1,2\}$, there is an approximation due to the
truncation of the screening functions outside $\tilde{\Omega}$,
and we have $\wt{\varphi}^{D\P,\w F}(\v
x_{\midx})\approx\varphi^{D\P,\w F}(\v x_{\midx})$. This
approximation is accurate as long as $\v x_{\midx}$ lies in
$\Omega$ and $\tilde{\Omega}$ is large enough, a statement that
will be made precise in the following sections.

\subsection{Outline of the method}
\label{sec:SE_outline}

So far, we have presented a unified formula for the
representation of the Fourier space sum with different
periodicities, i.e.~\eqref{eq:fourier_space_ewald}. Later, in
section~\ref{sec:disc}, we will introduce a uniform grid and
discretize the formulation. Before we do that, we want to give a
schematic outline of the SE method, without explicitly
introducing the grid. Of particular interest here is the way the
window function is introduced into what will become the gridding,
scaling and gathering steps, cf.~the numbered list at the
beginning of section~\ref{sec:SE} above. In
section~\ref{sec:windows}, we will review four relevant window
functions that are used in different Ewald methods.

Let $\W(\v x)$ be a window function with Fourier transform
$\wh{\W}(\v k)$ and consider the trivial identity
\begin{align}
  1 \equiv \wh{\W}(\v k) [\wh{\W}(\v k)]^{-2} \wh{\W}(\v k).
  \label{eq:window_split}
\end{align}
We will insert this into the Fourier space
sum~\eqref{eq:fourier_space_ewald} and use one factor $\wh{\W}$
for the gridding step, the second $\wh{\W}$ for the gathering
step, and $\wh{\W}^{-2}$ in the scaling step. Note that
\eqref{eq:fourier_space_ewald} can be rewritten as
\begin{equation}
  \wt{\varphi}^{D\P,\w F}(\v x_{\midx}) = 4\pi
  \mathcal{L}\left[
    \wh{\W}(\v k)
    e^{\ii \v k \cdot \v x_\midx}
    e^{-\lvert \v k \rvert^2/4\xi^2} [\wh{\W}(\v k)]^{-2} \wh{G}^{D\P}(\v k)
    \sum_{\nidx=1}^N q_{\nidx}
    \wh{\W}(\v k)
    e^{-\ii \v k \cdot \v x_\nidx}
  \right].
  \label{eq:fourier_space_ewald_rearranged}
\end{equation}
We start by defining
\begin{align}
  \wh{H}(\v k):=\sum_{\nidx=1}^N q_{\nidx}\wh{\W}(\v k)e^{-\ii \v
  k\cdot\v x_\nidx},
  \label{eq:fourier_gridding}
\end{align}
which is the Fourier transform of
\begin{align}
  H(\v x) = \sum_{\nidx=1}^N q_{\nidx}\W(\v x-\v x_\nidx)_\ast.
  \label{eq:gridding}
\end{align}
Here, $(\cdot)_{\ast}$ denotes that periodicity is implied in the
periodic directions. Eq.~\eqref{eq:gridding} will become the
gridding step in the final method, in which $H(\v x)$ is
evaluated on the uniform grid.

We furthermore define
\begin{align}
  \wh{\wt{H}}(\v k):= e^{-|\v k|^2/4\xi^2} [\wh{\W}(\v k)]^{-2} \wh{G}^{D\P}(\v k)
  \wh{H}(\v k),
  \label{eq:scale}
\end{align}
which will become the scaling step. Note that the scaling step
will depend on the selected window function. Also recall that
$\wh{G}^{D\P}$ here is the modified Green's function defined by
Eqs.~\eqref{eq:G_hat_3d}--\eqref{eq:G_hat_0d}.
Eq.~\eqref{eq:scale} allows us to write
\eqref{eq:fourier_space_ewald_rearranged} as
\begin{align}
  \wt{\varphi}^{D\P,\w F}(\v x)=4\pi\mathcal{L}\left[\wh{\W}(\v
  k)\wh{\wt{H}}(\v k)e^{\ii \vk\cdot \v v} e^{\ii \vka\cdot
  \v w}\right].
  \label{eq:int}
\end{align}
Applying the convolution theorem to \eqref{eq:int}, we obtain an
integral form for the Fourier-space potential evaluated at target
point $\v x_\midx$,
\begin{align}
  \wt{\varphi}^{D\P,\w F}(\v x_\midx) =
  4\pi\int_{\mathbb{R}^{3-D}}\int_{[0,L)^D}\wt{H}(\v v,\v w)\W(\v
  x_\midx-\v x)_\ast \, \w d\v v\w d \v w.
  \label{eq:complete_integral}
\end{align}
This form, once discretized, will become the gathering step. We
will discretize the integral in \eqref{eq:complete_integral}
using the trapezoidal rule. Note that if $\W$ is smooth and have
compact support, the integral can be computed with spectral
accuracy. Furthermore, note that in the case when target points
and source points coincide, the window function will be evaluated
in the same points in both \eqref{eq:gridding} and
\eqref{eq:complete_integral}, so its values can be cached.

In the above procedure, we obtain $\wh{H}$ from $H$ and $\wt{H}$
from $\wh{\wt{H}}$ via a Fourier transform and inverse Fourier
transform, respectively. In the 2d- and 1d-periodic cases these
transforms are \textit{mixed}: a discrete Fourier transform in
each periodic direction and an approximation to the continuous
Fourier integral in each non-periodic direction. We shall return
to this discussion later, in section~\ref{sec:aft}.

\subsection{Discretization}
\label{sec:disc}

Let us now introduce a uniform Cartesian grid on $\Omega = [0,L)^3$
with $M$ subintervals in each spatial direction ($M$ being an
even integer), and define the step size $h=L/M$. The charges of
the point sources are spread to the uniform grid using a suitable
window function in the gridding step (cf.~\eqref{eq:gridding}).
Let the support of the window function have length $Ph$ in each
spatial direction, where $P$ is a positive even integer such that
$P \leq M$ (i.e.\ $P$ is the number of grid subintervals within
the support of the window).

In the non-periodic directions, the box~$\Omega$ must be extended
to length $\tilde{L} = L + \delta_L$, as mentioned in
section~\ref{sec:vico}. The extended box must at least accomodate
the support of the window function, which means that $\delta_L
\geq Ph$. However, as shown in Ref.~\onlinecite{Klinteberg2017}
(figure~3), this may not be sufficient since the screening
function itself, cf.~\eqref{eq:gamma}, must also have decayed
sufficiently where the box ends. When the window function is a
Gaussian, an expression for $\delta_L$ that takes both the window
function and screening function into account can be derived
\cite{Klinteberg2017}, but for a general window function this may
not be possible. Instead, we will require $\delta_L \geq \lambda
P h$, where $\lambda \geq 1$ depends on $D$ and the window function and is
determined through experiments (the value of $\lambda$ is given in
section~\ref{sec:params}).
The uniform grid must also be extended to $\tilde{M} = M +
\delta_M$ in the non-periodic directions. To keep the grid
spacing the same in all directions, we must have
$\delta_L/\delta_M=h$. To ensure both this and that $\tilde{M}$
becomes an even integer, as well as that $\delta_L \geq \lambda P
h$, we first set $\tilde{M} = 2\lceil (M + \lambda P)/2 \rceil$
and then set $\tilde{L} = h \tilde{M}$. This defines the extended
box $\tilde{\Omega}=[-\tfrac{1}{2}\delta_L,L+\tfrac{1}{2}\delta_L)^3$,
where $\delta_L = \tilde{L} - L$.
In the $D=0,1,2$ cases, further upsampling is needed in the free
directions, which is the topic of the next section.

\subsection{Adaptive FFT and upsampling of Fourier integrals}
\label{sec:aft}

So far, we have treated the singular case ($\vk=\v0$ and the
free-space case) by introducing modified Green's functions
(section~\ref{sec:vico}), which remove the singularity. The
sampling rate in Fourier space is $1/\Delta k = L / 2\pi$ in the
periodic directions ($k = n \Delta k$ with $n = -M/2, \ldots,
M/2-1$) and would, after discretizing the Fourier integrals, in
the absence of upsampling be $1/\Delta \kappa = \tilde{L} / 2\pi$
in the free directions ($\kappa = n \Delta \kappa$ with $n =
-\tilde{M}/2, \ldots, \tilde{M}/2-1$).
However, the modified Green's functions are oscillatory and thus
requires a larger sampling rate than $\tilde{L} / 2\pi$ to resolve.
Furthermore, when $\vk$ is small but non-zero for $D=1,2$, the
Green's function $\wh{G}^{D\P}([\vk,\vka]) = 1/(\lvert \vk
\rvert^2 + \lvert \vka \rvert^2)$ varies rapidly for small
$\vka$, and thus also needs an increased sampling rate. One
option would be to apply a global upsampling factor $s_g > 1$ such
that $1/\Delta \kappa = s_g \tilde{L} / 2 \pi$ for all periodic
modes, and the grid would be of size $M^D (s_g \tilde{M})^{3-D}$.
However, for larger values of $\vk$, the integrand varies
more slowly and thus no upsampling is
needed.\cite{Shamshirgar2017} Based on this observation, we have
developed an \emph{adaptive Fourier transform}
(AFT) that applies different upsampling factors in the free
directions for different periodic modes. Upsampling is here
achieved by zero-padding in real space, which increases the
sampling rate in Fourier space.

\begin{figure}[htbp]
  \centering\includegraphics{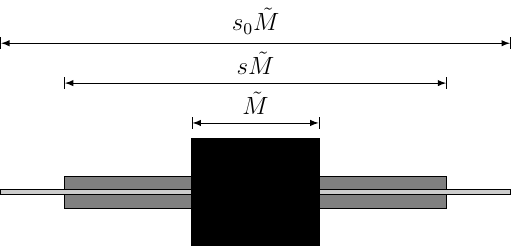}
  \caption{Schematic representation of the adaptive Fourier
  transform in two dimensions (vertical direction periodic and
  horizontal direction free).}
  \label{fig:upsampling}
\end{figure}

The AFT uses an upsampling factor $s_0$ for the $\vk=\v0$ mode
and free-space case, and an upsampling factor $s$ for non-zero
periodic modes satisfying $k_i \leq (2\pi/L) \nl$ for
$i=1,\ldots,D$ in the $D=1,2$ cases; higher modes are not
upsampled. Parameters to be determined are $s_0$, $s$ and $\nl$.
A schematic representation is given in
figure~\ref{fig:upsampling} in two dimensions. For $D=1,2$, the
AFT reduces the computational cost significantly compared to global
upsampling. Note that no upsampling is needed for $D=3$. For
$D=0$, the upsampling factor $s_0$ is applied globally; however,
part of the cost can be hidden in a precomputation step,
described in section~\ref{sec:precomputation}.

It has been shown that for free-space Poisson problems in $d=3-D$
dimensions, $s_0 \geq 1+\sqrt{d}$ is required to account for the
oscillatory behaviour of the modified Green's functions and
accurately compute the aperiodic convolution.\cite{Klinteberg2017}
In practice, we select $s_0 = 2$, $2.5$ and $2.8$ for $D=2$, $1$
and $0$, respectively.

The following is valid for $D=1,2$ only. To define which non-zero
periodic modes to upsample, let $\nl$ be a positive integer such
that $\nl\ll \kintinf := M/2$ and define the two sets
\begin{equation}
  \mathbb{I}:=\lbrace \vk \in \mathcal{K}^D \setminus \v 0:
                      |k_i| \leq \frac{2\pi}{L} \nl,\hspace{.5em} i = 1, \ldots, D \rbrace,
  \label{eq:I}
\end{equation}
and
\begin{equation}
  \mathbb{J}:=\lbrace \vk \in \mathcal{K}^D \setminus \mathbb{I}:
                      |k_i| \leq \frac{2\pi}{L} \kintinf,\hspace{.5em} i = 1, \ldots, D \rbrace.
  \label{eq:J}
\end{equation}
The set $\mathbb{I}$ contains periodic $\vk$-vectors in a box
around the zero vector (to be upsampled with factor $s$), while
$\mathbb{J}$ contains the vectors outside this box (which will
not be upsampled). The selection of $s$ and $\nl$ is treated
later in section~\ref{sec:params}. The AFT computes the Fourier
transform $\v w \rightarrow \vka$ in the free directions with
adaptive upsampling factor
\begin{equation}
  s_f(\vk) = \begin{cases}
    s_0, & \vk = \v 0, \\
    s, & \vk \in \mathbb{I}, \\
    1, & \vk \in \mathbb{J},
  \end{cases}
  \label{eq:upsamplings}
\end{equation}
i.e.\ $1/\Delta \kappa = s_f(\vk) \tilde{L}/2\pi$ and $\kappa = n
\Delta \kappa$ with $n = -s_f(\vk)\tilde{M}/2, \ldots,
s_f(\vk)\tilde{M}/2-1$. Recalling the
notation introduced in definition~\ref{def:definition_mixed}, we
now present the AFT/AIFT algorithms. These algorithms are the
generalized versions of those introduced for $D=1$ in
Ref.~\onlinecite{Shamshirgar2017}.

\begin{breakablealgorithm}
\caption{Adaptive Fourier transform - AFT for $D=1,2$}
\small
\label{alg:afft}
\begin{algorithmic}[1]
\INPUT Grid representation of sources $H$, upsampling factors $s_0, s$, and $\mathbb{I}$, $\mathbb{J}$ sets, cf.\ \eqref{eq:I} and \eqref{eq:J}.
\State Apply an FFT on $H$ in the periodic $\v v$ directions to compute $\wh{H}(\vk,\v w)$.
\State Pad $\wh{H}(\boldsymbol{0},\v w)$ with zeros in the free $\v w$ directions with upsampling factor $s_0$ and apply an FFT to compute $\wh{H}(\boldsymbol{0},\vka)$.
\State Pad $\wh{H}(\mathbb{I},\v w)$ with zeros in the free $\v w$ directions with upsampling factor $s$ and apply an FFT on $\vk\in\mathbb{I}$ to compute $\wh{H}(\mathbb{I},\vka)$.
\State Apply an FFT on $\wh{H}(\mathbb{J},\v w)$ with no upsampling to compute $\wh{H}(\mathbb{J},\vka)$.
\OUTPUT Adaptive Fourier transform of $H$, $\wh{H}$.
\end{algorithmic}
\end{breakablealgorithm}

Note that if $\nl=\kintinf$, all non-zero Fourier modes are
upsampled and therefore, steps 3 and 4 can be merged. If we also
have $s=s_0$, the AFT is the same as a plain three-dimensional
FFT.

Now, assume that as a result of applying the AFT algorithm,
$\wh{H}(\boldsymbol{0},\vka)$, $\wh{H}(\mathbb{I},\vka)$ and
$\wh{H}(\mathbb{J},\vka)$ are available. The adaptive inverse FFT
(AIFT) algorithm, which computes $H$ from $\wh{H}$, can easily be
found by reversing the steps of the AFT algorithm.

\begin{breakablealgorithm}
\caption{Adaptive inverse Fourier transform - AIFT for $D=1,2$}
\small
\label{alg:aifft}
\begin{algorithmic}[1]
\INPUT $\wh{H}(\boldsymbol{0},\vka)$, $\wh{H}(\mathbb{I},\vka)$ and $\wh{H}(\mathbb{J},\vka)$, grid size $M$ and $\mathbb{I},\mathbb{J}$ sets, cf.\ \eqref{eq:I} and \eqref{eq:J}.
\State Apply an IFFT on $\wh{H}(\mathbb{J},\vka)$ in the free $\v w$ directions to compute $\wh{H}(\mathbb{J},\v w)$.
\State Apply an IFFT on $\wh{H}(\mathbb{I},\vka)$ in the free $\v
  w$ directions and truncate the solution to $M$ grid points in each direction to compute $\wh{H}(\mathbb{I},\v w)$.
\State Apply an IFFT on $\wh{H}(\boldsymbol{0},\vka)$ in the free
  $\v w$ directions and truncate the solution to $M$ grid points
  in each direction to compute $\wh{H}(\boldsymbol{0},\v w)$.
\State Merge $\wh{H}(\mathbb{I},\v w)$, $\wh{H}(\mathbb{J},\v w)$ and $\wh{H}(\boldsymbol{0},\v w)$ to construct $\wh{H}(\vk,\v w)$.
\State Apply an IFFT on $\wh{H}(\vk,\v w)$ in the periodic $\v v$ directions to compute $H(\v v,\v w)$.
\OUTPUT Adaptive inverse Fourier transform of $\wh{H}$, $H$.
\end{algorithmic}
\end{breakablealgorithm}

Once again, if $\nl=\kintinf$ and $s=s_0$, all steps can be
merged into a three-dimensional inverse FFT followed by
truncating the result to $M^3$ grid points. It is also worth mentioning that for $D=3,0$, the AFT/AIFT is simply a three dimensional FFT/IFFT on a grid of size $M^3$ and $(s_0\tilde{M})^3$ respectively.

\subsection{Precomputation for the free-space case}
\label{sec:precomputation}

In the previous section we discussed how the AFT/AIFT algorithms
can be used to accelerate the computations compared to global
upsampling, for the $D=1,2$ cases. However, in the free-space
case this does not apply since there is no periodic directions
(and thus the AFT becomes equivalent to global upsampling).
Therefore, in the free-space case, all modes must be upsampled by
at least $s_0 \geq 1 + \sqrt{3} \approx 2.73$, applied to the
extended grid size $\tilde{M} = 2 \lceil (M+\lambda P)/2\rceil$
in all three directions. As stated before, upsampling is needed
to resolve the oscillatory modified Green's function
\eqref{eq:G_hat_0d} and to compute an aperiodic convolution.
However, the convolution only requires $s_0=2$; the rest is
needed to resolve the Green's function. By precomputing an
effective Green's function and truncating it in real space, we
can thus reduce zero-padding to a factor of 2 in each
direction.\cite{Klinteberg2017} The precomputation algorithm are
given in Algorithm~\ref{alg:precomp}. The output
$\wh{G}^{0\P}\subR$ can be stored and reused in the scaling step.
Note that it depends only on $\tilde{L}$ and $\tilde{M}$, not on
the locations or charges of the point sources.

\begin{breakablealgorithm}
\caption{Precomputation step for $D=0$}
\small
\label{alg:precomp}
\begin{algorithmic}[1]
  \INPUT Extended box side length $\tilde{L}$ and grid size
  $\tilde{M}$.
  \State Set $s_0=2.8$ and calculate $\wh{G}^{0\P}(k_1,k_2,k_3)$
  from \eqref{eq:G_hat_0d} for $k_i = 2 \pi n / (s_0 \tilde{L})$,
  \newline
  $n = -s_0\tilde{M}/2, \ldots, s_0\tilde{M}/2-1$, $i=1,2,3$.
  \State Apply a 3d IFFT to compute $G^{0\P}$ on a grid of size
    $(s_0 \tilde{M})^3$.
  \State Truncate $G^{0\P}$ to obtain $G^{0\P}\subR$ on a grid of
    size $(2 \tilde{M})^3$.
  \State Apply a 3d FFT to get $\wh{G}^{0\P}\subR$.
  \OUTPUT Effective Green's function $\wh{G}^{0\P}\subR(k_1, k_2,
  k_3)$ for $k_i = \pi n / \tilde{L}$, $n=-\tilde{M}, \ldots,
  \tilde{M}-1$.
\end{algorithmic}
\end{breakablealgorithm}

\subsection{Summary of the Spectral Ewald algorithm}
\label{sec:SE_summary}

In Algorithm~\ref{alg:spectral_ewald}, we finally present the
unified Spectral Ewald method to compute the approximate
Fourier-space part of the electrostatic potential, given by
\eqref{eq:fourier_space_ewald}, with arbitrary periodicity. Note
that in the free-space case ($D=0$), precomputation of
$\wh{G}^{0\P}\subR$ is done once according to
Algorithm~\ref{alg:precomp}; the result is then used in step~3 of
Algorithm~\ref{alg:spectral_ewald}.

\begin{breakablealgorithm}
\caption{Spectral Ewald method - $k$-space algorithm}
\small
\label{alg:spectral_ewald}
\begin{algorithmic}[1]\setcounter{ALG@line}{-1}
\INPUT Charge locations $\v x_{\nidx}\in[0,L)^3$ and charges
  $q_{\nidx}$, $\nidx=1,\ldots,N$, decomposition parameter~$\xi$,
  grid size~$M$, upsampling factors $s_0, s$, maximum upsampled
  Fourier mode $\nl$, window function support size $P$,
  number of periodic directions $D$.
  \State Set $h=L/M$, $\tilde{M}=2\lceil (M + \lambda P)/2 \rceil$,
  $\tilde{L}=h \tilde{M}$, $R=\sqrt{3-D}\tilde{L}$, $\delta_L =
  \tilde{L}-L$.
  \State (Gridding) Introduce a uniform grid on $[0,L)^D\times
  [-\tfrac{1}{2}\delta_L,L+\tfrac{1}{2}\delta_L)^{3-D}$ with $M^D \times \tilde{M}^{3-D}$ points.
  Evaluate $H([\v v,\v w])$ on this grid according to
  \eqref{eq:gridding}.
  \State (FFT) If $D=3$, apply an $M^D$-point FFT to compute
  $\wh{H}(\vk)$. If $D=2,1$, apply an AFT with parameters $s_0$,
  $s$ and $\nl$ to compute $\wh{H}([\vk, \vka])$. If $D=0$, apply
  an $(2\tilde{M})^{3-D}$-point FFT to compute $\wh{H}(\vka)$.
  \State (Scaling) Use \eqref{eq:scale} to obtain $\wh{\wt{H}}$.
  If $D=0$, use the precomputed $\wh{G}^{0\P}\subR$.
  \State (IFFT) Apply an IFFT (if $D=3,0$) or AIFT (if $D=2,1$) to
  compute $\wt{H}$.
  \State (Gathering) Compute the integral in \eqref{eq:complete_integral} at target points (same as source points) using the trapezoidal rule.
\OUTPUT Approximation to the Fourier space part of the potentials
  $\varphi^{D\P, \w F}(\v x_{\midx})$, $\midx=1,\ldots,N$.
\end{algorithmic}
\end{breakablealgorithm}

\section{Window functions}
\label{sec:windows}

Here, we review some of the most relevant window functions that
appear in electrostatic calculations and specifically in PME
methods. We also include a recent window function introduced by
Barnett \etal \cite{Barnett2019}. For a complete survey of the
classical window functions, the reader is directed to
Ref.~\onlinecite{Harris1978}.

The following window functions are presented in one dimension. In three dimensions, the corresponding window function $\W(\v x)$ can be obtained as a tensor product
\begin{align*}
  \W(\v x) = \W_0(x)\W_0(y)\W_0(z),
\end{align*}
where $\W_0$ is the one-dimensional window function.

\textbf{Gaussian window.}
Gaussians are the traditional window functions of the SE method%
\cite{Klinteberg2017,Lindbo2011a,Lindbo2012,Shamshirgar2017},
and they have several important properties. First and foremost,
their Fourier transforms (needed in the scaling step) are readily
available, and both the window and its Fourier transform are fast
to compute. Moreover, they are smooth and decay rapidly in
Fourier space. However, they do not have compact support and
therefore have to be truncated in practice. The truncated
Gaussian window function is defined as
\begin{equation}
  \W_\text{G}(x) = \begin{cases}
    e^{-\alpha (x/w)^2}, & |x| \leq w, \\
    0, & \text{otherwise},
  \end{cases}
  \label{eq:gaussian}
\end{equation}
where $\alpha>0$ is a shape parameter and $w > 0$ is the
half-width of the window function ($w = Ph/2$ with $P$ as in
section~\ref{sec:disc}). (In previous
work\cite{Klinteberg2017,Lindbo2011a,Lindbo2012,Shamshirgar2017}, the shape
parameter $m$ has been used, related to $\alpha$ through $\alpha
= m^2/2$.) The shape parameter $\alpha$ controls
the truncation level of the Gaussian, described further in
section~\ref{sec:approx}. Note that in higher dimensions, the
window is truncated outside of a cube.

The Fourier transform of the truncated Gaussian is given by
\begin{equation}
  \wh{\W}_\text{G}(k) = \sqrt{\frac{\pi}{\alpha}} w e^{-k^2w^2/4\alpha}
  \frac{\erf\!\left(\sqrt{\alpha} + \frac{\ii kw}{2\sqrt{\alpha}} \right)
  + \erf\!\left( \sqrt{\alpha} - \frac{\ii kw}{2\sqrt{\alpha}}
  \right)}{2},
\end{equation}
where $\erf$ is the error function. In practice, however, this
expression is never used; instead, we use the Fourier transform
of the untruncated Gaussian, i.e.
\begin{equation}
  \wh{\W}_\text{G,untrunc}(k) =
  \sqrt{\frac{\pi}{\alpha}} w e^{-k^2w^2/4\alpha},
\end{equation}
which is taken into account in the error analysis (see
section~\ref{sec:approx}).

\textbf{Cardinal B-spline window.} Another type of window
function that is used in FFT-based methods are cardinal
B-splines\cite{Essmann1995}. The B-spline of order 2 is defined
as
\begin{align*}
M_2(x)=
  \begin{cases}
  1-|x-1|, & 0\leq x\leq 2, \\
  0, & \text{otherwise},
  \end{cases}
\end{align*}
and for order $p>2$ is defined recursively as
\begin{align}
M_p(x) = \frac{x}{p-1}M_{p-1}(x)+\frac{p-x}{p-1}M_{p-1}(x-1).
\label{eq:bspline}
\end{align}
This window has finite support and is easy and fast to implement.
Moreover, its Fourier transform is available analytically.
B-splines have polynomial degree of smoothness and consequently,
if B-splines of low order are used in an FFT-based method, the
FFT grid size must be increased significantly to achieve high
accuracy.

\begin{figure}[htbp]
  \centering
  \includegraphics{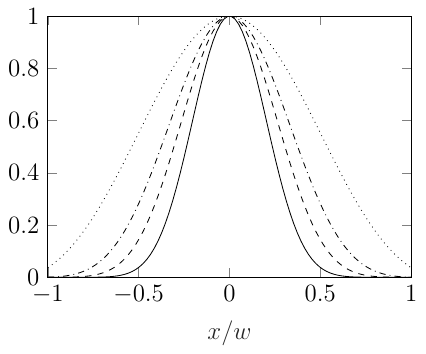}
  \caption{The KB window function with $\beta=5,10,15,25$, from top to bottom.}
  \label{fig:different_beta}
\end{figure}

\textbf{Kaiser-Bessel window.} The Kaiser-Bessel (KB) window
function\cite{Kaiser1980} is defined by
\begin{equation}
  \W_\text{KB}(x) = \begin{cases}
    \dfrac{I_0\left(\beta\sqrt{1-(\frac{x}{w})^2}\right)}{I_0(\beta)}, & |x|\leq w, \\
    0, & \text{otherwise},
  \end{cases}
  \label{eq:kaiser_window}
\end{equation}
where $I_0(\cdot)$ is the zeroth-order modified Bessel
function of the first kind and $\beta > 0$ is a shape parameter.
The KB window is shown with different $\beta$ in
figure~\ref{fig:different_beta}. We describe how $\beta$ should
be selected to minimize approximation errors in section~\ref{sec:approx}.
The Fourier transform of the window is available in closed form
as
\begin{equation}
  \wh{\W}_\text{KB}(k) = \dfrac{2w \,
  \sinh(\sqrt{\beta^2-k^2w^2})}{I_0(\beta)\sqrt{\beta^2-k^2w^2}}.
  \label{eq:kaiser_ft}
\end{equation}

The KB window is an approximation to a family of
prolate-spheroidal wave functions of order zero. It has been
shown that the prolate-spheroidal wave functions provide an
orthonormal basis which is optimal for the representation of
functions whose Fourier transforms are compactly supported
\cite{Osipov2013}. In addition, and to our interest in this
paper, they require a significantly smaller width~$w$ compared to
the Gaussian window to achieve the same target accuracy, thus
reducing the computational effort.
Potts \etal
\cite{Potts2003} used this window function in order to develop a
fast summation algorithm and recently Nestler \cite{Nestler2016}
showed that using the KB window function, the
resulting algorithm is more accurate than the method using
B-splines for homogeneous systems. More recently, Gao \etal
\cite{Gao2017} used the same window function in their simulation
and arrived at the same conclusion for non-homogeneous systems.
The main drawback of the KB window is that it is expensive
to compute. We address this by approximating it by a piecewise
polynomial, described in section~\ref{sec:polynomial-window}.

\textbf{Exponential of semicircle window.} The ``exponential of
semicircle'' (ES) window was recently introduced by Barnett \etal
\cite{Barnett2019} as an approximation to the KB function that
avoids evaluation of Bessel functions. The ES window is defined
by
\begin{equation}
  \W_\text{ES}(x) = \begin{cases}
   \dfrac{e^{\beta\sqrt{1-(\frac{x}{w})^2}}}{e^\beta}, & - w\leq x\leq w, \\
    0, & \text{otherwise}.
  \end{cases}
  \label{eq:expsemicirc_window}
\end{equation}
This window can achieve nearly as high precision as the KB window
with the same width, but is cheaper to evaluate. However, unlike
the other window functions introduced in this section, its
Fourier transform is not known analytically and thus has to be
computed numerically. Another drawback is that the ES window is
in fact not differentiable at the endpoints $x = \pm w$.

\begin{figure}[htbp]
  \centering
  \begin{subfigure}[b]{0.47\textwidth}
    \centering
    \includegraphics[width=\textwidth]{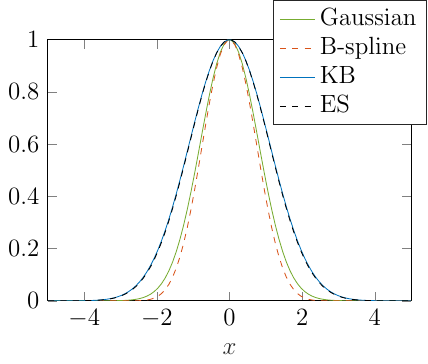}\vspace*{1.5mm}
  \end{subfigure}%
  \hfill%
  \begin{subfigure}[b]{0.49\textwidth}
    \centering
    \includegraphics[width=\textwidth]{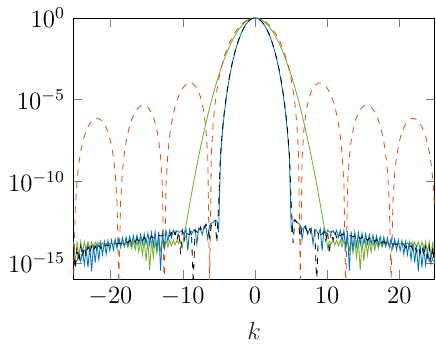}
  \end{subfigure}%
  \caption{(Left) An example of the Gaussian~\eqref{eq:gaussian},
  B-spline~\eqref{eq:bspline} of order $p=6$,
  KB~\eqref{eq:kaiser_window} and
  ES~\eqref{eq:expsemicirc_window} windows, all scaled to one.
  (Right) Decay of Fourier transforms of the window functions.
  For the Gaussian, KB and ES windows, $w=6$, $\alpha=28$ and
  $\beta=30$.}
  \label{fig:windows}
\end{figure}

In figure~\ref{fig:windows} we give an example of the Gaussian
window \eqref{eq:gaussian}, B-spline window \eqref{eq:bspline} of
order $p=6$, KB window \eqref{eq:kaiser_window} and ES window
\eqref{eq:expsemicirc_window} together with
their Fourier transforms, all scaled to one. For the KB and ES windows we choose
$w=6$ and $\beta=30$, and for the Gaussian window we set $w=6$
and $\alpha=28$. The parameters for the Gaussian, KB and ES
windows are selected such that they are truncated at the same
level approximately. The figure shows that the Gaussian decays
faster than the KB and ES windows in real space (left), but
slower in Fourier space (right). The KB and ES windows are very
similar in both real space and Fourier space.

The main advantage of the ES window is that it is a
cheaper-to-evaluate alternative to the KB window. However, using
the piecewise polynomial approximation described below, we can
evaluate the KB window itself to desired accuracy at low cost. We
will therefore have no need to consider the ES window, nor the
B-spline window, in the remainder of this paper.

\subsection{Polynomial approximation and the PKB window}
\label{sec:polynomial-window}

The idea is simply to approximate the window function by a
piecewise polynomial, with as high accuracy as needed to reach
the desired target accuracy in the SE method.
The inspiration for this comes from the FINUFFT
library\cite{Barnett2019}, where the idea is applied to the ES
window. The idea can be used with any (sufficiently smooth)
window function; here we apply it to the KB window. We
call the resulting piecewise polynomial window function the
``polynomial Kaiser-Bessel'' (PKB) window.
Note that since the function $I_0(\beta \sqrt{1-z^2})$, cf.\
\eqref{eq:kaiser_window}, is analytic in the whole complex plane,
it should lend itself well to polynomial approximation.

\begin{figure}[htbp]
  \centering
  \begin{subfigure}[b]{0.47\textwidth}
    \centering
    \includegraphics[width=\textwidth]{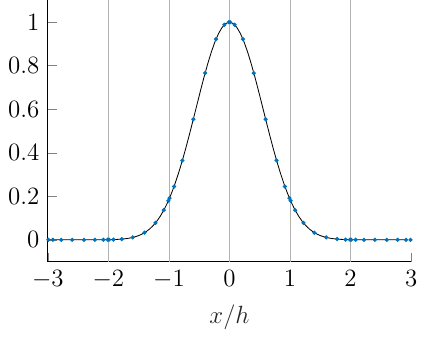}\vspace*{1.5mm}
  \end{subfigure}%
  \hfill%
  \begin{subfigure}[b]{0.49\textwidth}
    \centering
    \includegraphics[width=\textwidth]{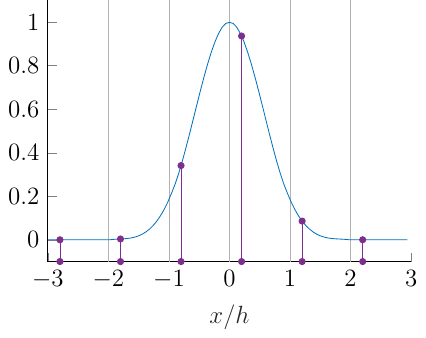}
  \end{subfigure}%
  \caption{(Left) In each subinterval, the window function is
  sampled in $\nu+1$ Chebyshev points (shown here for $P=6$
  subintervals and $\nu=7$) and interpolated by a polynomial of
  degree $\nu$. (Right) The piecewise polynomial is evaluated in
  $P$ grid points (purple), each with the same offset within its
  subinterval.}
  \label{fig:polynomial-construct}
\end{figure}

The support of the window function, which is of length $2w=Ph$
(cf.\ section~\ref{sec:disc}),
is divided into $P$ subintervals each of length $h$.
On each subinterval, the window function is interpolated by a
polynomial of degree $\nu$ (selection of $\nu$ is discussed
in section~\ref{sec:params}).
The interpolation is performed by first mapping each subinterval
$[x_i, x_i+h]$ to the interval $[-1, 1]$, i.e.\
\begin{equation}
  [x_i, x_i+h] \ni x \mapsto \bar{x} = 2
  \left(\frac{x-x_i}{h}\right) - 1 \in [-1,1].
\end{equation}
The window function in each subinterval is then sampled in
$\nu+1$ Chebyshev points
\begin{equation}
  \bar{x}_k = \cos \! \left( \frac{\pi (k - 1/2)}{\nu+1} \right), \qquad
  k=1,\ldots,\nu+1,
\end{equation}
illustrated in the left part of
figure~\ref{fig:polynomial-construct}, and subsequently the polynomial
\begin{equation}
  \label{eq:polynomial-ansatz}
  p_\nu(\bar{x}) = \sum_{j=0}^\nu c_j \bar{x}^j
\end{equation}
is fitted to the window function by solving the linear system
\begin{equation}
  \label{eq:polynomial-linsys}
  p_\nu(\bar{x}_k) = \W\!\left( \frac{\bar{x}_k+1}{2} h + x_i
  \right), \qquad k=1,\ldots,\nu+1
\end{equation}
for the $\nu+1$ coefficients $\{c_j\}$. The piecewise polynomial
is allowed to be discontinuous where two subintervals meet.

When the window function $\W(\v x)$ is to be evaluated in the gridding and
gathering steps, the window is centered at one of the point
sources and must be evaluated on the uniform grid. Thus,
the evaluation points all have the same offset within the
subintervals, shown in the right part of
figure~\ref{fig:polynomial-construct}. This fact makes the
evaluation efficient, since it means that $\bar{x}$ in
\eqref{eq:polynomial-ansatz} is in fact the same for all
subintervals; only the coefficients $\{c_j\}$ differ between
subintervals. The polynomial \eqref{eq:polynomial-ansatz} is
evaluated using Horner's rule. (In the scaling step, the
$\wh{\W}^{-2}$ factor is evaluated directly using
\eqref{eq:kaiser_ft}, i.e.\ no polynomial approximation is used.)

Since we use the monomial basis $\{\bar{x}^j\}$ in
\eqref{eq:polynomial-ansatz}, the matrix of the linear system
\eqref{eq:polynomial-linsys} becomes a $(\nu+1) \times (\nu+1)$
Vandermonde matrix evaluated in the Chebyshev points. The
condition number of this matrix is approximately $0.6 \times
10^{0.38 \nu}$, which seems acceptable considering that we will
see in section~\ref{sec:params} that $\nu=9$ (corresponding to a condition number of
$1.6 \times 10^3$) will be sufficient for full precision in the
SE method. In practice, the choice of polynomial basis does not
seem to influence the accuracy greatly, as long as Chebyshev
sampling points are used (for example, the Chebyshev basis gives
almost the same result as the monomial basis). However, note that
the choice of sampling points is important, and equidistant
sampling points would give a significantly more ill-conditioned
problem than Chebyshev points.

\begin{figure}[htbp]
  \centering
  \begin{subfigure}[b]{0.49\textwidth}
    \centering
    \includegraphics[width=\textwidth]{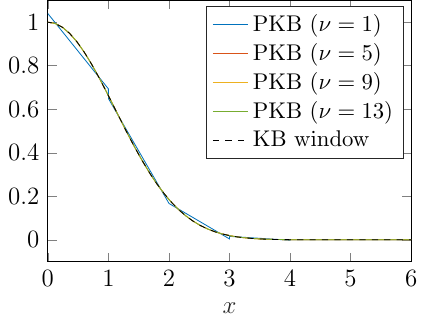}\vspace*{1.5mm}
  \end{subfigure}%
  \hfill%
  \begin{subfigure}[b]{0.49\textwidth}
    \centering
    \includegraphics[width=\textwidth]{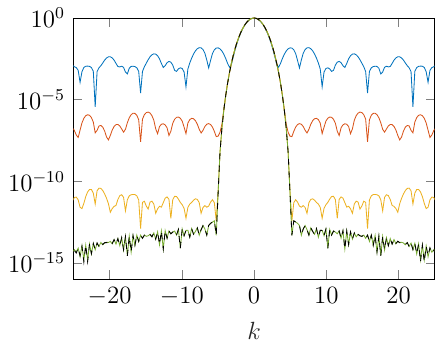}
  \end{subfigure}%
  \caption{(Left) PKB windows with different
  degrees~$\nu=1,5,9,13$ and (right) their Fourier transforms.
  The KB window which is approximated is shown as a dashed black
  curve. Here, $w=6$, $P=12$ and $\beta=30$. Only $x \geq 0$ is
  shown in the left part since all windows are even functions.}
  \label{fig:polynomial-window}
\end{figure}

The PKB window is determined by the number of subintervals~$P$,
shape parameter~$\beta$ (both given by the underlying KB window)
and degree~$\nu$. Examples of PKB windows with different
degrees~$\nu$ are given in figure~\ref{fig:polynomial-window}.
However, we show in section~\ref{sec:approx}
that $\beta$ can be related to $P$, and in
section~\ref{sec:params} that $\nu$ too can be related to $P$.
Thus, the PKB window will in the end be uniquely determined by $P$.

\section{Errors in the Spectral Ewald method}
\label{sec:errors}

\subsection{Truncation error}
\label{sec:truncation}
Truncation errors in Ewald methods are introduced due to the
truncation of interactions in the real space sum or truncation of
Fourier modes in the $k$-space sum. As a measure of accuracy, we
define the root mean square (rms) error as
\begin{equation}
  \erms =
  \left(\frac{1}{N}\sum_{\midx=1}^N|\wt{\varphi}_\midx-\varphi^\ast(\v
  x_\midx)|^2\right)^{1/2},
  \label{eq:rms_error}
\end{equation}
where $\varphi^\ast$ denotes the exact or well converged
potential. The error \eqref{eq:rms_error} may be evaluated in the
real-space part $\varphi^{D\P,\w R}$ or Fourier-space part
$\varphi^{D\P,\w F}$ separately. The magnitude of the truncation
errors can be perfectly estimated using error estimates by Kolafa and
Perram\cite{Kolafa1992}, which suggest that the error in the real
space sum is
\begin{equation}
  \erms^\text{R} \approx \sqrt{Q} \xi^{-2} (L r_\c)^{-3/2} e^{-(\xi r_\c)^2},
  \label{eq:KP_real_space}
\end{equation}
and the error in the $k$-space sum is
\begin{equation}
  \erms^\text{F} \approx \sqrt{Q} \xi \pi^{-2} \kintinf^{-3/2}
  e^{-[\pi \kintinf / (\xi L)]^2},
  \label{eq:KP_fourier_space}
\end{equation}
where $Q = \sum_{\nidx=1}^N q_{\nidx}^2$ and $\kintinf = M/2$.
Therefore, setting an absolute error tolerance $\erms$ and an
Ewald decomposition parameter~$\xi$, the cut-off radius should be
selected as
\begin{equation}
  r_\c = \frac{\sqrt{3}}{2\xi} \sqrt{W \left( \frac{4}{3}
  \left[ \frac{Q}{\xi L^3 \erms^2} \right]^{2/3} \right)},
  \label{eq:KP_rc}
\end{equation}
and the maximum wavenumber as
\begin{equation}
  \kintinf = \frac{\sqrt{3}\xi L}{2\pi} \sqrt{W \left(
  \frac{4}{3}
  \left[
    \frac{Q}{\pi\xi L^3\erms^2}
  \right]^{2/3}
  \right)},
  \label{eq:KP_kinf}
\end{equation}
where $W$ is the Lambert $W$ function (defined as the solution to
$W(x) e^{W(x)} = x$).

Assume that $\erms$ and $\xi$ are fixed and that the density
$N/L^3$ and charge density $Q/L^3$ are constant while $N$
increases. Then $r_\c$ as given by \eqref{eq:KP_rc} is fixed,
while $\kintinf \propto L$. Truncating the real-space
interactions beyond $r_\c$ reduces the computational complexity
for the real-space sum from $\order{N^2}$ to $\order{N}$,
with a constant that depends on the number of particles within a
ball of radius $r_\c$. For the $k$-space sum, the gridding and
gathering steps clearly have complexity $\order{N}$, while the
FFTs have complexity $\order{M^3 \log(M^3)}$ and the scaling step
$\order{M^3}$. Since $M^3 \propto \kintinf^3 \propto L^3
\propto N$, the $k$-space sum, and thus the whole algorithm,
scales as $\order{N \log(N)}$.

\subsection{Approximation error}
\label{sec:approx}

Approximation errors arise due to (a) the evaluation of \eqref{eq:complete_integral} with the trapezoidal rule using truncated window functions and (b) approximating Fourier integrals. 
We have noted earlier that the quadrature error of the Fourier integrals can be controlled 
by upsampling the grid using upsampling factors $s$ and $s_0$ and
the parameter $\nl$. In section~\ref{sec:params} we will discuss how to select these parameters. 
In this section, we focus only on approximation errors due to (a).

In Ref.~\onlinecite{Lindbo2011a}, where Gaussians are used as
window functions, the authors derive the approximation error
estimate
\begin{align}
  C_\text{G} \left(e^{-(\pi/2)^2P^2/\alpha} + \erfc(\sqrt{\alpha})\right),
  \label{eq:approx_error_Gaussian}
\end{align}
where $\alpha$ is the shape parameter of the Gaussian window and
$P$ is the size (i.e.\ number of subintervals) of the support of
the window in each direction, cf.~Theorem 3.1 in
Ref.~\onlinecite{Lindbo2011a} (where $m=\sqrt{2\alpha}$).
The constant $C_\text{G}$ may depend on the solution
$\varphi^{D\P,\w F}$, but not on $P$ or $\alpha$.
In \eqref{eq:approx_error_Gaussian}, the first term estimates the
quadrature error and the second term estimates the window
function truncation error.
Balancing both terms and using the approximation $\erfc(\sqrt{x})
\approx e^{-x}$,
one obtains $\alpha = (\pi/2) P c^2$, where $c^2$ is a heuristically
inserted constant which in practice is set slightly below unity
(in this paper we use $c^2=0.91$).

Also for the KB window function, the shape parameter~$\beta$
strongly influences the accuracy of the algorithm (but not the
runtime for a fixed $P$). Selection of the shape parameter for
the KB window has been discussed in different contexts for
instance in Refs.~\onlinecite{Kaiser1980}, \onlinecite{Potts2003}
and \onlinecite{Gao2017}. There is no universally optimal shape
parameter; rather, the optimal choice depends on the context in
which the window function is used. Here, we use numerical
experiments to determine a near-optimal shape parameter for use
in the SE method, in the sense that it nearly minimizes the
approximation error for all $P$. We heuristically find the
approximation error estimate
\begin{align}
  C_\text{KB} \left(e^{-2\pi P^2/\beta} + \erfc(\sqrt{\beta})\right)
  \label{eq:approx_error_KB}
\end{align}
for the KB window, where the constant $C_\text{KB}$ may depend on
the solution but not on $P$ or $\beta$. Note that this estimate
is very similar to the one provided in Eq.~\eqref{eq:approx_error_Gaussian}, and again
the first term estimates the quadrature error while the second
estimates the window truncation error. Using the approximation
$\erfc(\sqrt{x})\approx e^{-x}$, we find that the shape parameter
that balances both terms is given by $\beta = \sqrt{2\pi} P
\approx 2.5 P$. Note that, as for the Gaussian window, the shape
parameter is only a function of $P$.

\begin{figure}[htbp]
  \centering
  \begin{subfigure}[b]{0.49\textwidth}
    \centering
    \includegraphics[width=\textwidth]{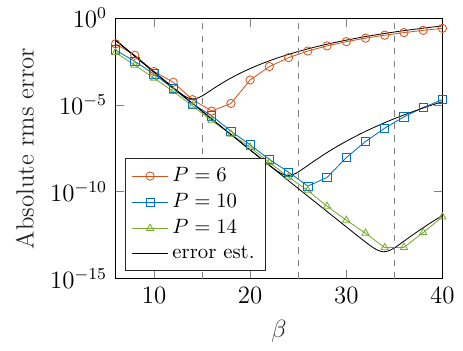}
  \end{subfigure}\hfill%
  \begin{subfigure}[b]{0.49\textwidth}
    \centering
    \includegraphics[width=\textwidth]{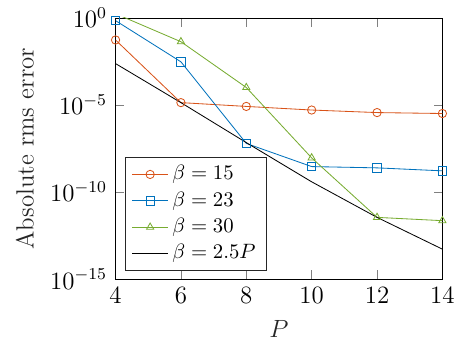}
  \end{subfigure}
  \caption{(Left) The rms error in
  evaluating~\eqref{eq:fourier_space_ewald} in the 3d-periodic
  case as a function of $\beta$, and the error
  estimate~\eqref{eq:approx_error_KB} with $C_\text{KB}=100$. Vertical dashed
  lines show $\beta=2.5P$. (Right) The rms error in
  evaluating~\eqref{eq:fourier_space_ewald} as a function of
  the KB window support~$P$ for different shape
  parameters~$\beta$. In both figures, $N=100$ ($Q=31.4$), $L=1$,
  $M=28$ and $\xi=6.5$.}
  \label{fig:shape_parameter}
\end{figure}

In figure \ref{fig:shape_parameter} (left), we plot the rms error
in evaluating the triply periodic electrostatic potential
\eqref{eq:fourier_space_ewald} as a function of $\beta$ for
$P=6,10,14$, together with the
error estimate  \eqref{eq:approx_error_KB}. This figure
demonstrates that $\beta=2.5P$ coincides with the minimum
of the error curves for each $P$. Figure
\ref{fig:shape_parameter} (right) shows the rms error as a
function of $P$ for different shape parameters as well as the
near-optimal choice $\beta=2.5P$.

\begin{figure}[htbp]
  \centering
  \begin{subfigure}[b]{0.49\textwidth}
    \centering
    \includegraphics[width=\textwidth]{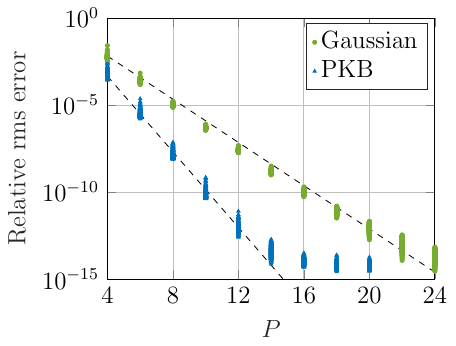}
  \end{subfigure}\hfill%
  \begin{subfigure}[b]{0.49\textwidth}
    \centering
    \includegraphics[width=\textwidth]{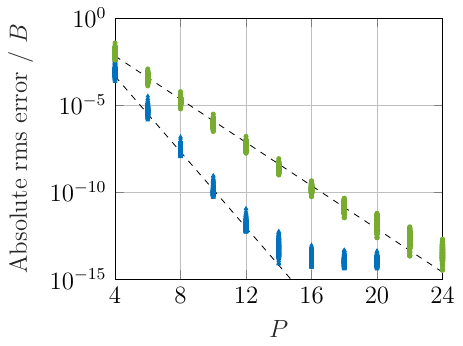}
  \end{subfigure}
  \caption{(Left) The relative rms error and (right) absolute rms
  error divided by $B = \sqrt{Q}f_c(\xi L)/L$, when evaluating
  \eqref{eq:fourier_space_ewald}, as a function of $P$. 480
  different uniformly randomly distributed systems are shown,
  given by all combinations of $N=100, 200, 400, 800$, $L=1, 5,
  10, 20, 40, 80$, $\xi L = 5, 15, 25, 30, 35$ and $D = 0, 1, 2,
  3$, for both the Gaussian and PKB window. $Q$ varies between 28 and 280.
  Other parameters are chosen such that other errors are
  negligible. Dashed lines show the relative error estimates $2
  e^{-(\pi/2) Pc}$ and $10 e^{-2.5 P}$ for the Gaussian and PKB
  window, respectively.}
  \label{fig:AB_test}
\end{figure}

Using the optimal shape parameters $\alpha = (\pi/2) P c^2$ and
$\beta = 2.5 P$, the error estimates
\eqref{eq:approx_error_Gaussian} and \eqref{eq:approx_error_KB}
can be simplified to $2 C_\text{G} e^{-(\pi/2) P c}$ and $2
C_\text{KB} e^{-2.5 P}$, respectively. Thus the approximation error is
controlled solely by $P$.
Moreover, in Ref.~\onlinecite{Shamshirgar2017}, we showed that
$C_\text{G} \approx A := \sqrt{Q \xi L}/L$, where
$Q = \sum_{\nidx=1}^N q_{\nidx}^2$. Here, we show that this can be
refined to $C_\text{G} \approx B$ and that $C_\text{KB} \approx 5
B$, where $B := \sqrt{Q} f_c(\xi L) / L$ with
\begin{equation}
  f_c(x) = e^{-c_1/x^2} (c_2 + c_3 x + c_4 x^2),
  \label{eq:fc-expression}
\end{equation}
where $c_1=12.62$, $c_2=0.8909$, $c_3=0.01411$ and $c_4=4.315
\times 10^{-5}$ were obtained by curve fitting. The right part of
figure~\ref{fig:AB_test} shows that the approximation error is
well estimated by $2 B e^{-(\pi/2) P c}$ and $10 B e^{-2.5 P}$
for the Gaussian and PKB window, respectively.
As shown by the left part of the figure, $B$ approximates the
magnitude of the potential (i.e.\ the relative error is
approximately the absolute error divided by $B$). Thus, scaling
the error estimates by $B$, one obtains estimates for relative
errors.

\begin{figure}[htbp]
  \centering
  \includegraphics{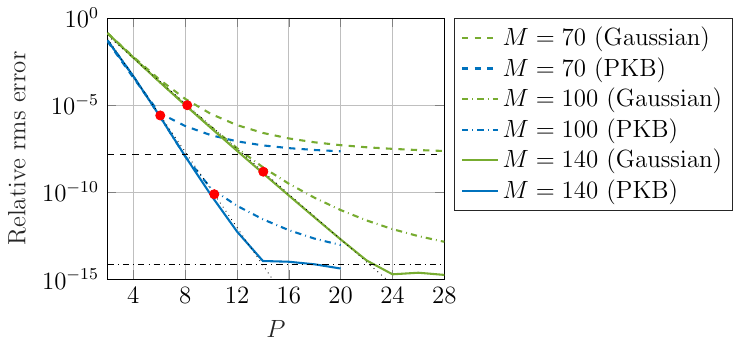}
  \caption{Relative rms error when evaluating
  \eqref{eq:fourier_space_ewald} in the triply periodic case, as
  a function of $P$. Parameters are $N=1000$, $L=5$, $\xi=6$;
  $Q=334$. The error is shown for the Gaussian and PKB windows for $M=70$
  (dashed curves), $M=100$ (dot-dashed curves) and $M=140$ (solid
  curves). According to \eqref{eq:KP_fourier_space}, these values of $M$
  correspond to relative truncation errors of $1.5\times
  10^{-8}$, $7.4 \times 10^{-15}$ and $1.7 \times 10^{-26}$,
  respectively, indicated by black horizontal lines (the last of
  these is not visible). Note that
  the error curves deviate from the approximation error estimates
  $2Be^{-(\pi/2)Pc}$ and $10Be^{-2.5P}$ (dotted black lines)
  significantly above the truncation error level. The red dots
  indicate where $P=f_\W(M/(\xi L))$, with $f_\W$ as in
  \eqref{eq:error_pollution}. The red dots for $M=140$ are below
  $10^{-15}$ and thus not visible.}
  \label{fig:error_pollution}
\end{figure}

We have observed that when $P$ and $M$ are both selected from
error estimates, i.e.\ when $2 B e^{-(\pi/2) P c}
\approx \erms \approx \eqref{eq:KP_fourier_space}$ for the
Gaussian window or $10 B e^{-2.5 P} \approx \erms \approx
\eqref{eq:KP_fourier_space}$ for the PKB window, the truncation
error and approximation error may pollute each other so that the
total error becomes larger than expected. We have heuristically
determined that this occurs when $P > f_{\W}(M/(\xi L))$, where the
function $f_{\W}$ depends on the window function and is
\begin{equation}
  f_{\W}(x) = \begin{cases}
    x^2 + 0.2 x + 2.25 & \text{for the Gaussian window}, \\
    0.7 x^2 + 0.2 x + 1.8 & \text{for the PKB window}.
  \end{cases}
  \label{eq:error_pollution}
\end{equation}
This is illustrated in figure~\ref{fig:error_pollution}.
In this situation, the error still decreases down to the
truncation error level when $P$ is increased, but
slower than predicted by the error estimates. We suggest
increasing $M$ by 5~\% and increasing the window support to $P+4$
when $P > f_{\W}(M/(\xi L))$, which should be enough to reach the
error expected from the estimates.

\section{Parameter selection}
\label{sec:params}

There are several parameters present already in the standard SE
method for the triply periodic case, and we have here added some
more parameters related to upsampling in the AFT and the
polynomial window function. Here, we provide a systematic
approach for selecting these parameters.

Let $\xi$ and an error tolerance $\erms$ be given. (The parameter
$\xi$ is in principle free, and is in practice used to balance
the runtime of the real-space and Fourier-space computations.)
Select $r_\c$ and $\kintinf$ using the Kolafa \& Perram
estimates \eqref{eq:KP_rc} and \eqref{eq:KP_kinf}, respectively.
This also sets the grid size in the periodic directions as
$M=2 \kintinf$.

\begin{figure}[htbp]
  \centering
  \includegraphics{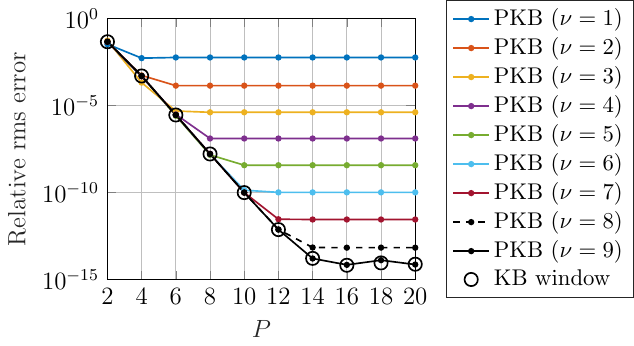}
  \caption{Relative rms error as a function of $P$ using the PKB
  window with different polynomial degrees~$\nu$ in the 3d-periodic
  case. Error using the KB window shown as circles. Parameters
  are $N=10^5$, $L=1$, $\xi=8$, $M=32$, and $\beta=2.5P$.}
  \label{fig:poldeg-vs-P}
\end{figure}

Either the Gaussian window \eqref{eq:gaussian} or KB window
\eqref{eq:kaiser_window} is selected, the latter evaluated
using the PKB window (section~\ref{sec:polynomial-window}).
To set $P$, the error estimates $\erms \approx 2 B e^{-(\pi/2) P c}$
and $\erms \approx 10B e^{-2.5 P}$ are used for the Gaussian and (P)KB window functions,
respectively; $P$ is rounded to an even integer. However, to
avoid error pollution as explained in section~\ref{sec:approx}, if $P >
f_{\W}(M/(\xi L))$ where $f_{\W}$ is given by
\eqref{eq:error_pollution}, then the values of $P$ and $M$ (and
thus $\kintinf$) determined above are modified by increasing $M$
by 5~\% and adding~4 to $P$ (unless $M$ and $P$ were set directly
and not through the error estimates).
In any case, the window is truncated at $w=Ph/2$ where $h=L/M$. Based on the
discussion in the previous section, the shape parameters are selected as
$\alpha = (\pi/2)Pc^2$ and $\beta=2.5P$ for the Gaussian and
(P)KB windows, respectively. For the Gaussian window, we use the
fast Gaussian gridding approach, see
Ref.~\onlinecite{Lindbo2011a} and references therein.
For the PKB window, the polynomial
degree~$\nu$ must be selected. As shown in
figure~\ref{fig:poldeg-vs-P}, $\nu=9$ is sufficient to reach the
same precision as the KB window for all $P$, and the rule
\begin{equation}
  \nu = \min \! \left(\frac{P}{2} + 2, ~ 9 \right)
\end{equation}
can be used to select $\nu$.

The remaining parameters only appear when $D=2,1,0$. In the free
directions, the grid is extended to $\tilde{M} = 2\lceil (M +
\lambda P)/2 \rceil$ and the box is extended to $\tilde{L} = h
\tilde{M}$, as discussed in section~\ref{sec:disc}.
The parameter $\lambda$ depends on the
periodicity and window function and is given in
table~\ref{tab:lambda}.
\begin{table}[htbp]
  \centering
  \caption{Value of parameter $\lambda$.}
  \label{tab:lambda}
  \begin{tabular}{l@{\quad}c@{\quad}c@{\quad}c@{\quad}c}
    \hline
    & $D=0$ & $D=1$ & $D=2$ & $D=3$ \\
    \hline
    Gaussian window & 1 & 1.5 & 1.5 & --- \\
    (P)KB window & 1.3 & 2.4 & 2.4 & --- \\
    \hline
  \end{tabular}
\end{table}
The reason that we select $\lambda$ larger for
the (P)KB window compared to the Gaussian window is that the
product $\lambda P$ then becomes approximately the same for both
windows for a given tolerance, so that the FFT grids (i.e.\
$\tilde{M}$ and $\tilde{L}$) and the upsampling parameters below
become similar for both windows.
Note that $\lambda$ does not necessarily have to be larger for
$D=1,2$ than for $D=0$, but increasing $\lambda$ slightly in the
former cases leads to a lower runtime since $s$ and $\nl$ can be
selected smaller, as seen below.

Based on the requirement $s_0 \geq 1+\sqrt{d}$, where $d$ is the
number of free directions, for the upsampling factor in the
free-space case and for the zero modes in the 2d- and 1d-periodic
cases, we select $s_0$ equal to 2, 2.5 and 2.8 in the
2d-periodic, 1d-periodic and free-space cases, respectively.
Finally, for the 1d- and 2d-periodic cases, we also need to
select $s$ and $\nl$ for the $\mathbb{I}$ set of the AFT,
cf.~section~\ref{sec:aft}. Based on a generalized form of an
estimate derived in Ref.~\onlinecite{Shamshirgar2017}, we set
(see Appendix~\ref{app:aft_params})
\begin{equation}
  s = \frac{M}{\tilde{M}} \left( 1 + \frac{1}{2\pi} \log \! \left(
  \frac{B}{2 \erms} \right) \right)
  \label{eq:estimate-s}
\end{equation}
and
\begin{equation}
  \nl = \left\lceil \frac{M}{\tilde{M}-M} \frac{1}{2\pi}
  \log \! \left( \frac{B}{2\erms} \right) - 1 \right\rceil.
  \label{eq:estimate-nl}
\end{equation}
In practice, $s_0$ and $s$ are adjusted such that $s_0 \tilde{M}$
and $s \tilde{M}$ are both even integers.

\section{Numerical results}
\label{sec:results}

Here, we provide numerical results of three kinds: (i) a
comparison between the new PKB window function and the classical
Gaussian window in the SE method, (ii) force computation with the
SE method, and (iii) a comparison with the FMM.
In all cases, we
measure the relative rms error, either in the full electrostatic
potential~\eqref{eq:potential_split} or in the Fourier-space
part~\eqref{eq:fourier_space_ewald}. All experiments are
performed on a single core of a machine with
an Intel Core i7-8700 CPU which runs at 4.6~GHz with 32~GB of memory.
Unless otherwise stated, parameters are selected according to
section~\ref{sec:params}.

The unified Spectral Ewald code 
will soon be available online\cite{se_github}. 
The core routines are written in C and dynamically linked and
called through a \matlab\ MEX interface. Fourier transforms are
computed using \matlab's \texttt{fft}, which is based on the FFTW
library \cite{Frigo1998}. The Spectral Ewald package was built
with the GNU C compiler, version
5.4.0. 
The implementation also allows for simulation of systems with non-cubic box shapes.

\subsection{Comparison between the Gaussian and PKB window functions}

In the first numerical experiment, we compare the Gaussian and
PKB window functions in terms of cost and accuracy, when
computing the Fourier-space part \eqref{eq:fourier_space_ewald}
of the potential. We choose a system of $N=10^5$ random
particles, with positions $x_\nidx$ uniformly distributed in a
unit box $\Omega = [0,1)^3$ and charges $q_\nidx$ uniformly
distributed in the interval $[-1, 1]$, with $Q=3.32 \times 10^4$.
Furthermore, we set $\xi=8$ and $M=32$ (corresponding to a
truncation error around $1.65 \times 10^{-17}$,
cf.~\eqref{eq:KP_fourier_space}).
The window support $P$ is varied, while other
parameters, such as upsampling factors, are selected as in
section~\ref{sec:params} with tolerance $\erms=1.5 \times
10^{-17}$. Thus, the window function approximation errors
(cf.~section~\ref{sec:approx}) should dominate, while other
errors should be around machine precision.

\begin{figure}[htbp]
  \centering
  \includegraphics[width=0.44\textwidth]{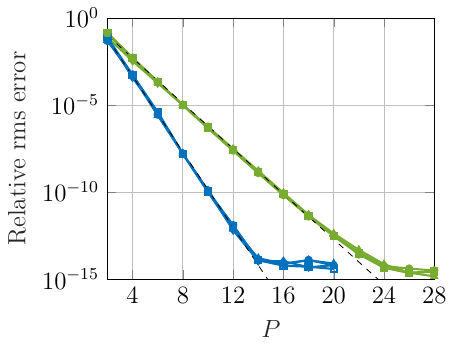}\hspace{5mm}%
  \includegraphics[width=0.45\textwidth]{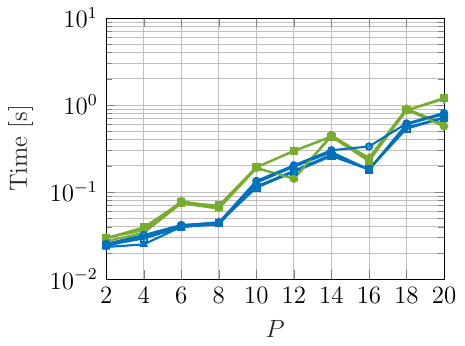}\\[1mm]
  \includegraphics{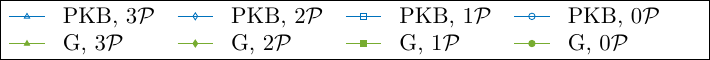}
  \caption{(Left) Relative rms error and (right) runtime for
  gridding and gathering as a function of $P$ using the PKB and
  Gaussian windows for the 3d-, 2d-, 1d- and 0d-periodic
  (free-space) cases. Parameters are $N=10^5$, $L=1$, $\xi=8$ and
  $M=32$.
  The other parameters are chosen such that other errors are
  negligible. Note that the runtime includes only the cost of the
  gridding and gathering steps, see
  Algorithm~\ref{alg:spectral_ewald}. The dashed lines in the
  left part are the error estimates $2 e^{-(\pi/2) P c}$ for the Gaussian and
  $10 e^{-2.5 P}$ for the PKB window.}
  \label{fig:efficiency}
\end{figure}

In figure~\ref{fig:efficiency} (left), the relative rms error is
plotted as a function of $P$. The result suggests that $n$ digits
of accuracy can be achieved with $P \approx n$ for the PKB
window (and machine precision is reached with $P \approx 14$),
while the Gaussian window requires $P \approx 1.6 n$ to achieve
$n$~digits (and reaches machine precision with $P \approx 24$).
Thus, we have approximately that $P_\mathrm{G} \approx 1.6P_\mathrm{PKB}$.

The runtime as a function of $P$ is demonstrated in
figure~\ref{fig:efficiency} (right) for both window functions.
The same system as in the left figure is used, and we only
include the runtime of the gridding and gathering steps.
This figure suggests that the PKB window is slightly faster than
the Gaussian window for a given $P$ in most cases. Since the PKB
window can also use a smaller $P$ than the Gaussian, it will
require less time to reach a given error.

\begin{figure}[htbp]
  \centering
  \includegraphics[width=0.45\textwidth]{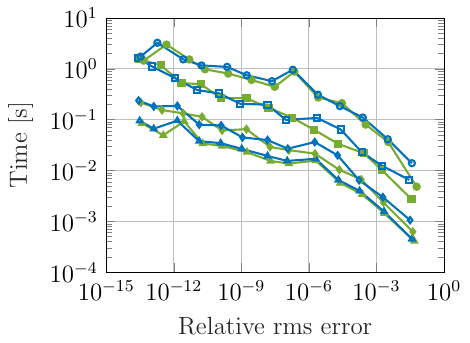}\hspace{5mm}%
  \includegraphics[width=0.45\textwidth]{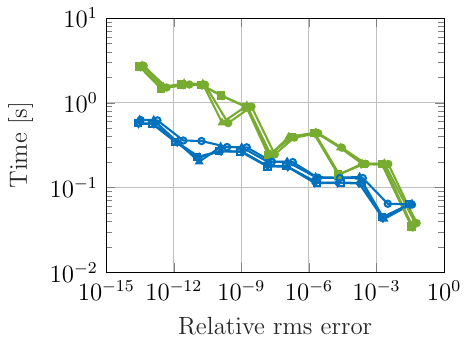}\\[1mm]
  \includegraphics{fig/window_legend}
  \caption{Runtime of (left) the FFT and scaling steps and
  (right) the gridding and gathering steps versus relative rms
  error, using the PKB and Gaussian windows with $D=0,1,2,3$
  periodic directions. Parameters are $N=10^5$, $L=10$ and
  $\xi=3$. Other parameters such as $M$ and $P$ are determined
  from error estimates as described in section~\ref{sec:params}.
  The free-space precomputation (section~\ref{sec:precomputation})
  is not included in the runtimes.}
  \label{fig:time_fftgrid_rms_kaiser}
\end{figure}

To study the effect of periodicity on the total computational
cost of the SE method, let us as a second experiment consider a
uniformly distributed system of $N=10^5$ particles in a box of
size $L=10$ (with the same charge distribution as in the previous
experiment, i.e.\ $Q=3.32 \times 10^4$). We select the
decomposition parameter $\xi=3$ and choose all other parameters
from the absolute error tolerance $\erms = 10^{-12}, 10^{-11}, \ldots,
10^{-1}, 10^0$ as described in section~\ref{sec:params}. In
figure~\ref{fig:time_fftgrid_rms_kaiser} the Fourier-space
runtime, excluding the free-space precomputation, is illustrated
separately for the FFT+scaling (left) and gridding+gathering
(right) steps.
The runtimes in the left plot are mainly functions of the grid size
$M$ and periodicity, while the results in the right plot are
functions of $N$ and $P$ (and hence the selected window function).

Evidently, for cases with the same type of periodic boundary
conditions, the runtime curves of the FFT+scaling steps are
essentially independent of the window function. (Naturally, this is due to
the fact that, as noted in section~\ref{sec:params}, the size of
the FFT grid is kept approximately the same for both window
functions.) Moreover, as the number of non-periodic directions
increases, the total number of grid points and consequently the cost of the
FFT+scaling steps grows. Nonetheless, the AFT algorithm reduces
the runtime significantly compared to full upsampling, such that
the 2d-periodic case has a similar (around a factor two larger)
computational cost to the 3d-periodic case. In the free-space
case, where the AFT algorithm cannot be used, we see the full
effect of upsampling in all three non-periodic directions
(although by a factor 2 rather than 2.8 due to the precomputation
step).
The right plot in figure~\ref{fig:time_fftgrid_rms_kaiser}
demonstrates that the PKB window function reduces the runtime of
the gridding+gathering steps by a factor 2--4 in most cases. This
is because $P$ can be selected smaller for the PKB window
compared to the Gaussian.

\begin{figure}[tbp]
  \centering
  \includegraphics[width=0.5\textwidth]{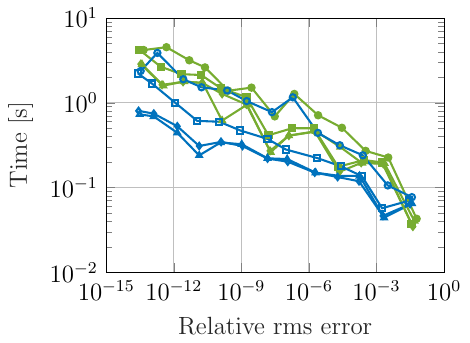}
  \caption{Total runtime (excluding precomputation) versus
  relative rms error, using the PKB and Gaussian windows with
  $D=0,1,2,3$ periodic directions. Parameters and legends are the
  same as in figure~\ref{fig:time_fftgrid_rms_kaiser}. The
  free-space precomputation (section~\ref{sec:precomputation}) is
  not included in the runtimes.}
  \label{fig:time_rms_kaiser}
\end{figure}

Finally, we sum the runtimes in both parts of the figure to
obtain the total runtime, as shown in
figure~\ref{fig:time_rms_kaiser}. It is clear that, compared to
the Gaussian window, the PKB window function reduces the total
runtime significantly in the 3d- and 2d-periodic cases, where the
gridding+gathering steps are the main bottleneck, while the
improvement is smaller in the free-space case where the
FFT+scaling steps dominate the cost. Thus, relatively speaking,
the total runtime is more affected by periodicity for the PKB
window than for the Gaussian window.

In Ref.~\onlinecite{Shamshirgar2017} we showed for the Gaussian
window that the total runtime of the 1d-periodic case is only
somewhat larger than the 3d-periodic case, due to the
employment of AFTs. The relative difference between the
1d-periodic and 3d-periodic cases is a bit larger for the PKB
window, due to the observation in the previous paragraph. With
the PKB window, the 1d-periodic case is at most three times more
costly than the 3d-periodic case, for very strict error tolerances.
In contrast, the runtime of the 2d-periodic case is very close to
that of the 3d-periodic case, for both window functions. For the
free-space case, the precomputation step reduces the cost of the
FFT+scaling steps. Even so, the free-space case is around four
times more expensive than the 3d-periodic case for strict
tolerances.

The precomputation step for the free-space case takes about as
much time as all other steps (gridding+gathering+FFT+scaling)
combined. The reason we exclude the precomputation from the total
runtime is that, in a time-dependent simulation where the system
is simulated for a long time, the precomputation is a one-time
cost and will therefore be negligible compared to the other steps
of the algorithm which must be carried out in every time step.

\subsection{Force computation}
\label{sec:force}

In molecular dynamics simulations, one is usually interested in
evaluating the force and energy along with the potential. The
Fourier-space part of the energy can be obtained simply by using
\begin{equation*}
  E^\text{F} = \sum_{\midx=1}^Nq_\midx \varphi^{D\P,\w F}(\v
  x_\midx).
\end{equation*}
To compute the Fourier-space force $\v F^\text{F}(\v x_\midx) = -
q_\midx \nabla_{\v x_\midx}
\varphi^{D\P,\w F}(\v x_\midx)$, the potential has to be differentiated with respect
to $\v x_\midx$. This can be done in two ways, both preserving the
spectral accuracy.
\begin{enumerate}
  \item
    By differentiating \eqref{eq:complete_integral} with respect
    to $\v x_\midx$, i.e.\
    \begin{equation}
      \nabla_{\v x_\midx} \wt{\varphi}^{D\P,\w F}(\v x_\midx) =
      4\pi\int_{\mathbb{R}^{3-D}}\int_{[0,L)^D}\wt{H}(\v v,\v w)
      [\nabla \W(\v x_\midx-\v x)_\ast] \, \w d\v v\w d \v w,
      \label{eq:force-method-1}
    \end{equation}
    Algorithm~\ref{alg:spectral_ewald} can be used
    to calculate the force, modifying only step~5 (gathering).
    This requires the derivative of the window function, and
    since the output from \eqref{eq:force-method-1} is a
    three-dimensional vector, the integration must be performed
    three times\cite{Shamshirgar2017}.
  \item
    Alternatively, the differentiation can be carried out in Fourier
    space, which amounts to multiplying the scaling step
    \eqref{eq:scale} by a factor $\ii \v k$. This modifies the
    output from step~3 in Algorithm~\ref{alg:spectral_ewald} to
    be a vector, so both step~4 and step~5 must be carried out
    three times.
\end{enumerate}
While method~2 is very easy to implement, it requires three times
as many IFFTs as method~1, and therefore method~2 is generally
avoided in molecular dynamics simulations. Furthermore, method~2
may require $M$ to be increased since the extra factor $\ii \v k$
makes the quantity $\wh{\wt{H}}(\v k)$, cf.~\eqref{eq:scale},
decay slightly slower in Fourier space. Method~1, on the other
hand, requires the gradient of the window function, and since the
window function is
a tensor product $\W(\v x) = \W_0(x) \W_0(y) \W_0(z)$, this means
that $\W_0'(x)$ must be available. When the PKB window is used,
we do not differentiate the polynomial approximation. Instead, we
construct a separate polynomial approximation of
$\W_\text{KB}'(x)$, cf.~\eqref{eq:kaiser_window}, using the same
procedure described in section~\ref{sec:polynomial-window} and
same values for $P$, $\beta$ and $\nu$.

\begin{figure}[bp]
  \centering
  \includegraphics[width=0.45\textwidth]{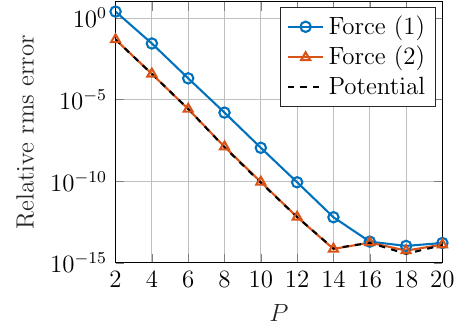}\hspace{8mm}%
  \includegraphics[width=0.45\textwidth]{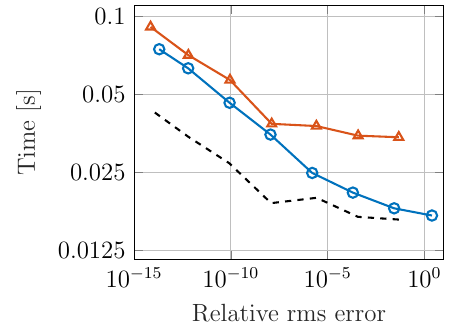}%
  \caption{Relative rms error versus (left) $P$ and
  (right) total runtime, when computing the Fourier-space force using methods~1
  and 2. Data for the potential is also shown for reference
  (dashed curve). All simulations are in the 3d-periodic case
  with the PKB window, with parameters $N=10^4$, $L=1$, $\xi=8$
  and $M=64$.}
  \label{fig:force}
\end{figure}

Here, we demonstrate both methods to compute the force. We
consider a system of $N=10^4$ uniformly distributed particles,
generated as before ($Q=3.35 \times 10^3$), in a box of size
$L=1$ with periodicity in all directions, i.e.\ $D=3$. We set
$\xi=8$ and $M=64$, and compute the forces acting on all
particles using methods~1 and 2, with the PKB window function.
The result is shown in figure~\ref{fig:force}.
The left plot shows that method~2 has the
exact same error curve as the potential (since $M$ is here large
enough for truncation errors to be negligible in both cases), while the error curve of
method~1 is offset by a constant factor. This reflects the fact
that method~1 uses a different window function in the gathering
step (namely the derivative of the original window function),
which requires a larger $P$ to resolve at the same error level as
the original window function. Here, $P$ needs to be increased by
2 in method~1 to reach the same error as method~2. (In principle,
it would be enough to increase $P$ in the gathering step, leaving
the gridding step unchanged. However, increasing $P$ in both the
gridding and gathering steps allows the structure of the scaling
step, with the factor $\wh{\W}^{-2}$, to be the same,
cf.~\eqref{eq:scale}.)

Despite the fact that $P$ needs to be slightly increased in
method~1, it is still faster than method~2 for all tolerances
(even more so for less strict tolerances), as
shown in the right plot of figure~\ref{fig:force}. The reason is
that method~2 requires three IFFTs rather than one (while both
methods require three integrations in the gathering step, unlike
the potential computation which only requires one). However, the
relative difference in runtime between methods~1 and 2 depends on
how dominant the cost of IFFTs are in the total algorithm, which
in turn depends on $N$, $M$ and $P$ (as well as the periodicity).

\subsection{Runtime comparison with FMM3D}
\label{sec:fmm3d}

In this section, we strive to demonstrate that the SE method is
competitive with other state-of-the-art fast summation methods.
To that effect, we compare the runtimes of the SE method and the FMM3D
library\cite{fmm3d}, a recent improved implementation of the Fast
Multipole Method (FMM)\cite{Cheng1999}. While both methods can be
applied to problems with different periodicities, the FMM is
fastest in the free-space case, whereas the SE method is fastest
in the triply periodic case. Therefore, to give each algorithm
optimal conditions we compare the SE method applied to a
triply-periodic problem with the FMM applied to a free-space
problem. The full potential \eqref{eq:potential} is computed in
both cases.

While the FMM is widely used and has been celebrated as one of
the top ten algorithms of the 20th century\cite{Dongarra2000},
it is less often used in molecular dynamics simulations e.g.\ in
GROMACS. The reason is that the FMM gives rise to discontinuities
in the potential and force which break the conservation of energy
in the system, although this violation of energy conservation is
at the selected error level.\cite{Shamshirgar2019} As shown in
Ref.~\onlinecite{Shamshirgar2019}, the FMM can be regularized to
alleviate this problem, albeit at an extra cost. However, we will
not consider this problem here.

It should be noted that the FMM3D library comes in two versions,
one ``easy to install'' version and one high-performance
optimized version, said to be up to three times faster than the
``easy'' version on some CPUs. Since our own Spectral Ewald code
is not highly optimized, we here use the ``easy'' version of
FMM3D. Furthermore, all comparisons are performed on a single
core, on the same machine as above.

\begin{figure}[htbp]
  \centering
  \includegraphics[width=0.35\textwidth]{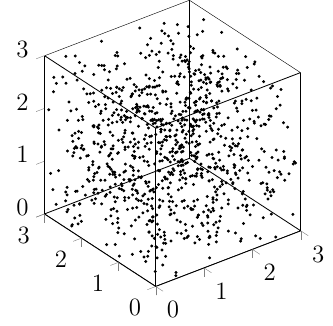}\hspace{1.5cm}%
  \includegraphics[width=0.35\textwidth]{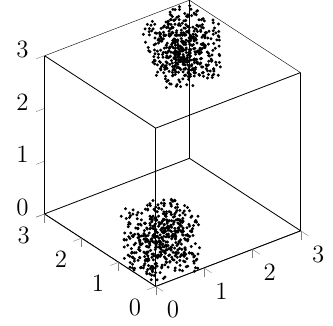}\\[1mm]
  \caption{Particle systems used in the comparison between the SE
  method and FMM3D: (Left) uniformly distributed particles and
  (right) isolated dense clouds. For the SE method, both problems
  are periodic in all three directions, while for FMM3D, both
  problems are free in all directions.}
  \label{fig:fmm_systems}
\end{figure}

We consider two different types of particle systems, each with
$N$~particles in a box of size $L=3$. In the first system (left
part of figure~\ref{fig:fmm_systems}), the particles are
uniformly distributed in the box. In the other (right part of
figure~\ref{fig:fmm_systems}), all particles are concentrated in
two corners of the box, thus forming isolated dense clouds
surrounded by empty space. Both the SE method and the FMM are
expected to be faster for the former, uniform, system than the
dense-cloud system, but the FMM is expected to handle the latter
system better than the SE method due to adaptivity. In both
cases, the charges are uniformly randomly selected from the
interval $[-1,1]$.

\begin{figure}[htbp]
  \centering
  \includegraphics[width=0.45\textwidth]{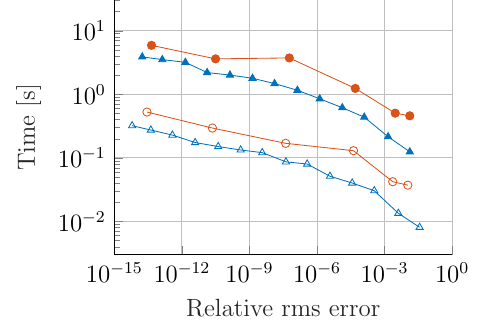}\hspace{5mm}%
  \includegraphics[width=0.45\textwidth]{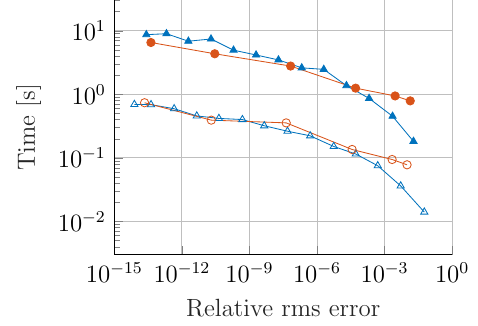}\\[1mm]
  \includegraphics{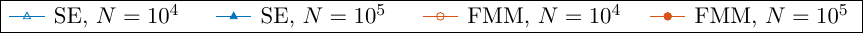}
  \caption{Total runtime versus relative rms error for the SE
  method and FMM3D applied to (left) the uniform system and
  (right) the dense-cloud system, both with $L=3$. For FMM3D, all
  problems are in free-space. For the SE method, all problems are
  triply periodic, the PKB window is used, $\xi$ is as stated in the
  text and other parameters selected from the error tolerance.}
  \label{fig:fmm_results}
\end{figure}

We consider both particle systems with $N=10^4$ ($Q=3.35 \times 10^3$) and
$N=10^5$ ($Q=3.32 \times 10^4$) particles.
For FMM3D the problem is considered in free-space. For the SE
method, the triply-periodic problem is considered, we use the
PKB window and select $\xi$ such that the runtimes of the
real-space and Fourier-space parts are balanced;
for the uniform system we get $\xi=10$ for
$N=10^4$ and $\xi=17$ for $N=10^5$, while for the dense-cloud
system we get $\xi=11$ for $N=10^4$ and $\xi=27$ for $N=10^5$.
For both methods, the absolute error tolerance is varied between
$10^{-14}$ and $10^0$. The total runtime versus relative rms
error for the full potential is shown in
figure~\ref{fig:fmm_results}. For the uniform
system, the SE method is slightly faster than FMM3D for all $N$
and tolerances considered here. For the dense-cloud system, the
methods are comparable for $N=10^4$, while FMM3D is slightly
faster than SE for $N=10^5$, as expected. However, for very loose
error tolerances the SE method seems to be faster also in this
case.

\section{Conclusions}
\label{sec:conclusions}

In this work, we have presented a unified approach to compute
the Fourier-space part of the Ewald sum for the electrostatic
potential, with arbitrary periodicity (triply, doubly, singly
periodic as well as free-space). We use the Spectral Ewald
method together with the recent idea of Vico \etal
\cite{Vico2016} to unify the treatment of all Fourier modes and
to utilize FFTs to accelerate the calculation. This approach has
already been used for the 1d-periodic \cite{Shamshirgar2017} and
free-space cases \cite{Klinteberg2017}, but not previously for the
2d-periodic case. In an attempt to compute the zero modes in this
case, Lindbo and Tornberg \cite{Lindbo2012} used
Chebyshev interpolation and showed that the extra cost of
calculating these modes is small compared to the total cost of
the SE method. In this paper we unified these modes into the
treatment of the other modes and thus completed the framework.

We also extended the idea of the adaptive Fourier transform,
first developed in Ref.~\onlinecite{Shamshirgar2017}, to the
2d-periodic case. With this, upsampling can be applied to a
small fraction of the Fourier modes, thus further reducing the
computational cost.

We compared the Gaussian window function, previously used in the
SE method, with a piecewise polynomial approximation of the
Kaiser-Bessel window function, here referred to as the PKB
window. The PKB window accelerates the gridding and gathering
steps of the method by reducing the number of points in the
support compared to the Gaussian window. With the PKB window
function, the support can be reduced by about 40~\% compared to
the Gaussian window for the same accuracy. The new window
function maintains the spectral accuracy observed with Gaussians,
and the method as a whole scales as $\order{N\log(N)}$ for $N$
sources. Implementing the new window function in the SE method,
we compared the resulting method with the existing algorithm that
has been developed in a series of papers
\cite{Lindbo2011a,Lindbo2012,Shamshirgar2017,Klinteberg2017}.

To compute the optimal shape parameter of the PKB window, we
numerically estimated the approximation error due to the
truncation of the Kaiser-Bessel window function and the
employment of the trapezoidal rule. We then showed that the
optimal shape parameter, and therefore the approximation error,
can be controlled using a single parameter, namely the number of
points in the support of the window function. The same parameter
controls the degree of the polynomial approximation.

Our numerical experiments show that the AFT algorithm can
reduce the runtime remarkably such that the cost of computing
doubly periodic and triply periodic cases are very similar, while
the singly periodic case is at most three times more expensive,
the most for very strict error tolerances. In order to accurately
compute Fourier integrals in the free-space case, the grid has to
be upsampled up to two times in each of the three directions.
Our results indicate that, when precomputing Fourier coefficients
of the modified Green's function, the free-space case is around
four times as expensive as the triply periodic
case for strict tolerances. We presented and compared two methods
to compute forces using the SE method, and also made a runtime comparison between
the SE method and the Fast Multipole Method to indicate that the
SE method is competitive.

The unified Spectral Ewald package for solving electrostatic
problems with different boundary conditions is accelerated using
OpenMP and vector intrinsics and will soon be available
online\cite{se_github}. 

\begin{acknowledgments}
  This work has been supported by the  G\"{o}ran Gustafsson
  Foundation for Research in Natural Sciences and Medicine and by
  the Swedish e-Science Research Center (SeRC). The authors
  gratefully acknowledge this support.
\end{acknowledgments}

\appendix

\section{Modified Green's function for doubly periodic case}
\label{app:G_hat_2d}

For completeness, we here provide an outline of the derivation of
Eq.~\eqref{eq:G_hat_2d}. Consider the one-dimensional Poisson
problem
\begin{equation}
  - \frac{\mathrm{d}^2\varphi}{\mathrm{d}z^2} (z) = f(z), \qquad z \in \mathbb{R},
  \label{eq:1d_Poisson}
\end{equation}
with boundary conditions $\varphi(z) \to 0$ as $\lvert z \rvert
\to \infty$. We select the symmetric Green's function $G(z) =
-\lvert z \rvert /2$ to \eqref{eq:1d_Poisson}. The Fourier
transform of the truncated Green's function is then
(cf.~section~\ref{sec:vico})
\begin{equation}
  \wh{G}\subR(\kappa) = \int_{-\infty}^\infty G(z)
  \, \w{rect}\!\left(\frac{z}{R}\right) e^{-\ii \kappa z} \, \w d
  z
  = \frac{1 - \cos(R|\kappa|) - R |\kappa| \sin(R|\kappa|)}{|\kappa|^2},
\end{equation}
and $\lim_{\kappa \to 0} \wh{G}\subR(\kappa) = -R^2/2$.

\section{Estimates for local upsampling in the adaptive Fourier transform}
\label{app:aft_params}

We here provide estimates for the local upsampling factor $s$ and
the size $\nl$ of the upsampling region, cf.~Eqs.~\eqref{eq:I}
and \eqref{eq:upsamplings} in section~\ref{sec:aft}. We start
from a generalization of Theorem~2 in
Ref.~\onlinecite{Shamshirgar2017}, which we state without proof:

\textit{Let $k, a > 0$ and $b \in \mathbb{R}$ and define
\begin{equation}
  F(k, a, b) = \int_{\mathbb{R}} f(\kappa; k, a, b) \, \w{d}
  \kappa, \qquad
  f(\kappa; k, a, b) = \frac{e^{-a(\kappa^2+k^2)}}{\kappa^2+k^2}
  e^{\ii b \kappa}.
  \label{eq:appB_integral}
\end{equation}
For any $\Delta \kappa > 0$, define the trapezoidal rule
approximation
\begin{equation}
  T_{\Delta \kappa}(k,a,b) = \Delta \kappa \sum_{j=-\infty}^\infty f(j\Delta \kappa; k, a, b).
\end{equation}
Then we have the error estimate
\begin{equation}
  \lvert T_{\Delta \kappa}(k, a, b) -
  F(k,a,b) \rvert \approx \frac{2\pi}{k}
  \frac{1}{e^{2\pi k / \Delta \kappa} - 1} \cosh(b k)
  =: H_{\Delta \kappa}(k,b).
\end{equation}}

We point out that $k$ may be thought of as the periodic
directions and $\kappa$ as the free directions; comparing
\eqref{eq:appB_integral} with
\eqref{eq:fourier_space_explicit_2P}, it is clear that this can
be directly applied to the doubly periodic case, with
$k^2=k_1^2+k_2^2$, $\kappa = \kappa_3$, $a = 1/4\xi^2$ and
$b=z_\midx-z_\nidx$. Noting that $\Delta \kappa = 2\pi / (s_f
\tilde{L})$, cf.~section~\ref{sec:aft}, and that $\lvert b \rvert
\leq L$, we can write
\begin{equation}
  H_{\Delta \kappa}(k,b) \leq L \bar{k}^{-1}
  \frac{1}{e^{2\pi (\tilde{L}/L) s_f \bar{k}} - 1} \cosh(2\pi \bar{k}),
\end{equation}
where $\bar{k} = Lk/(2\pi)$. Using the approximation $\cosh(x)
\approx \frac{1}{2} e^x$ and assuming that $e^{2\pi
(\tilde{L}/L) s_f \bar{k}} \gg 1$, we get
\begin{equation}
  H_{\Delta \kappa}(k,b) \lesssim \frac{L \bar{k}^{-1}}{2}
  e^{-2\pi \bar{k} [(\tilde{L}/L) s_f - 1]} =: \tilde{H}(\bar{k},
  s_f).
\end{equation}
Based on Eq.~\eqref{eq:upsamplings}, section~\ref{sec:aft}, we let
\begin{equation}
  s_f = \begin{cases}
    s_0, & \bar{k}=0, \\
    s, & 1 \leq \bar{k} \leq \nl, \\
    1, & \bar{k} \geq \nl+1.
  \end{cases}
\end{equation}
For $1 \leq \bar{k} \leq \nl$, the value of $\tilde{H}(\bar{k}, s)$
is maximal when $\bar{k}=1$. Demanding that $\tilde{H}(1, s) \leq
\varepsilon L$ for some $\varepsilon > 0$, we find that
\begin{equation}
  s \geq \frac{L}{\tilde{L}} \left( 1 - \frac{1}{2\pi} \log (2
  \varepsilon) \right).
  \label{eq:proto-estimate-s}
\end{equation}
To determine $\nl$, we furthermore demand that $\tilde{H}(\nl+1, 1) \leq
\varepsilon L/(\nl+1)$, which yields
\begin{equation}
  \nl \geq - \frac{L}{\tilde{L}-L} \frac{1}{2\pi} \log(2\varepsilon) - 1.
  \label{eq:proto-estimate-nl}
\end{equation}
Noting that $\varepsilon$ is a relative error tolerance here, we
get an absolute error tolerance $\erms$ by setting $\varepsilon =
\erms/B$, with $B$ as above Eq.~\eqref{eq:fc-expression}. Also
noting that $L/\tilde{L} = M/\tilde{M}$,
Eqs.~\eqref{eq:proto-estimate-s} and \eqref{eq:proto-estimate-nl}
lead to \eqref{eq:estimate-s} and \eqref{eq:estimate-nl},
respectively. While derived for the double periodic case, these
estimates are empirically found to work well also in the singly
periodic case.

\bibliographystyle{aipnum4-1}
\bibliography{library}

\end{document}